\documentclass[AMA,Times1COL]{WileyNJDv5} %

\hypersetup{colorlinks=true}

\usepackage[english]{babel}
\usepackage{amsmath}
\usepackage{amssymb}
\usepackage{amsthm}
\usepackage{nicefrac}
\usepackage[dvipsnames]{xcolor}
\usepackage{microtype}
\usepackage{comment}
\usepackage{graphicx} %
\usepackage{subfigure}
\usepackage{matlab-prettifier}
\usepackage{tikz}
\usepackage{pgfplots}
\pgfplotsset{compat=1.18}

\usepackage{subcaption}%
\usepackage{listings}%
\usepackage{matlab-prettifier}%
\definecolor{MatlabCellColour}{RGB}{252,251,220}%

\definecolor{mycolor1}{rgb}{0.00000,0.44700,0.74100}%
\definecolor{mycolor2}{rgb}{0.85000,0.32500,0.09800}%
\definecolor{mycolor3}{rgb}{0.92900,0.69400,0.12500}%
\definecolor{mycolor4}{rgb}{0.49400,0.18400,0.55600}%

\theoremstyle{definition}
\newtheorem{example}{Example}

\articletype{Research article}%

\received{Date Month Year}
\revised{Date Month Year}
\accepted{Date Month Year}
\journal{Journal}
\volume{00}
\copyyear{2023}
\startpage{1}

\definecolor{MyDarkGreen}{rgb}{0.0,0.4,0.0}
\begin{document}

\title{A Riemannian Optimization Approach for Finding the Nearest Reversible Markov Chain}

\author[1]{Fabio Durastante}

\author[1]{Miryam Gnazzo}

\author[1]{Beatrice Meini}

\authormark{\textsc{Durastante}, \textsc{Gnazzo}, \textsc{Meini}} %
\titlemark{A Riemannian Optimization Approach for Finding the Nearest Reversible Markov Chain}

\address[1]{\orgdiv{Mathematics Department}, \orgname{University of Pisa}, \orgaddress{\state{PI}, \country{Italy}}}

\corres{Corresponding author Miryam Gnazzo \email{miryam.gnazzo@dm.unipi.it}}

\presentaddress{Largo Bruno Pontecorvo, 5 – 56127 Pisa, PI. Italy.}

\fundingInfo{Italian Ministry of University and Research (MUR) through the PRIN 2022 ``Low-rank Structures and Numerical Methods in Matrix and Tensor Computations and their Application'' code: 20227PCCKZ MUR D.D. financing decree n. 104 of February 2nd, 2022 (CUP I53D23002280006) and through the PRIN 2022 ``MOLE: Manifold constrained Optimization and LEarning'',  code: 2022ZK5ME7 MUR D.D. financing decree n. 20428 of November 6th, 2024 (CUP B53C24006410006). }

\abstract[Abstract]{We address the algorithmic problem of determining the reversible Markov chain $\widetilde X$ that is closest to a given Markov chain $X$, with an identical stationary distribution. More specifically, $\widetilde X$ is the reversible Markov chain with the closest transition matrix, in the Frobenius norm, to the transition matrix of $X$. To compute the transition matrix of $\widetilde X$, we propose a novel approach based on Riemannian optimization. Our method introduces a modified multinomial manifold endowed with a prescribed stationary vector, while also satisfying the detailed balance conditions, all within the framework of the Fisher metric. We evaluate the performance of the proposed approach in comparison with an existing quadratic programming method and demonstrate its effectiveness through a series of synthetic experiments, as well as in the construction of a reversible Markov chain from transition count data obtained via direct estimation from a stochastic differential equation.}

\keywords{Markov Chain, Reversible, Riemannian Optimization, Stochastic matrices}

\jnlcitation{\cname{%
\author{Durastante F.}, 
\author{Gnazzo M.}, and 
\author{Meini B.}. 
\ctitle{A Riemannian Optimization Approach for Finding the Nearest Reversible Markov Chain} \cjournal{\it Numer. Linear Algebra Appl.} \cvol{2025;00(00):1--18}.}}

\maketitle

\renewcommand\thefootnote{\fnsymbol{footnote}}
\setcounter{footnote}{1}

\section{Introduction}

An ergodic Markov chain, with transition probability matrix $P$ and stationary vector $\boldsymbol{\pi}$, is said to be reversible if $D_{\boldsymbol{\pi}} P = P^\top D_{\boldsymbol{\pi}}$,  where $D_{\boldsymbol{\pi}}$ is the diagonal matrix with the components of $\boldsymbol{\pi}$ on its diagonal. Reversibility is a fundamental property in Markov chain theory, with profound consequences for both theoretical insights and practical implementations. Many algorithms---particularly those in statistical physics\cite{Seneta}, computational biology\cite{PhysRevE.73.046126}, and Markov chain Monte Carlo (MCMC) methods\cite{MR1397966}--- exploit the reversibility property to secure desirable convergence and stability characteristics. For instance, in MCMC one constructs a reversible Markov chain whose stationary distribution coincides with a target law, thereby enabling the design of optimal sampling schemes for approximating complex distributions\cite{Wu2015184}. In computational biology, the assumption of reversibility permits the faithful reconstruction of molecular dynamics from partial observations, yielding robust models of system kinetics\cite{Fackeldey201773}. More generally, reversibility underpins the detection of almost-invariant aggregates in nearly uncoupled chains—systems typified by extended metastable excursions interrupted by rare transitions—and justifies lumping techniques for dimensionality reduction\cite{Tifenbach2013120,Grone2008}.

In practice, however, the transition matrix of a Markov chain is often derived from noisy measurements or estimated from data, and hence may not satisfy the reversibility condition\cite{PhysRevE.93.033307}. This discrepancy motivates the need to approximate a given (possibly non-reversible) Markov chain with a reversible one that is, in an appropriate sense, closest to the original.

In some cases, reversibility is enforced by employing ``reversibilizations'' that utilize closed-form modifications based on transition acceptance/rejection probabilities\cite{Choi_Wolfer,Choi}. Although such methods yield a reversible chain, they do not provide a quantitative measure of how close the resulting chain is to the original. 

An alternative strategy involves formulating and solving a quadratic programming (QP) problem\cite{Nielsen2015483} to compute the reversible Markov chain whose transition matrix is the nearest—in the Frobenius norm—to a given stochastic matrix. 

Our approach follows this spirit and exploits Riemannian optimization techniques to fully employ the geometric structure inherent in the set of transition matrices of reversible Markov chains with a prescribed stationary distribution. In particular, our construction relies on the introduction of a new Riemannian manifold, $\mathcal{M}_{\boldsymbol{\pi}}$, which generalizes the manifold of stochastic matrices with fixed stationary distribution\cite{DurastanteMeini} for approximating the stochastic $p$th root of a stochastic matrix, which itself generalizes the multinomial manifolds\cite{Douik8861409}. 
Given the stationary distribution $\boldsymbol{\pi}$ and the transition probability matrix $A$ of a Markov chain,
the optimization problem that we solve is the computation of
the closest matrix in the manifold  $\mathcal{M}_{\boldsymbol{\pi}}$—in the Frobenius norm—to the matrix $A$.
Since it is possible to construct a second-order geometry for this manifold, we can employ efficient second-order optimization methods for the optimization task.
Moreover, the proposed algorithm automatically detects any transient states and restricts its analysis to the remaining recurrent subchain, thus enhancing both efficiency and robustness.
In cases where the underlying Markov chain decomposes into several ergodic classes, our approach 
processes each one independently. By isolating and solving smaller, decoupled subchains, we achieve further substantial computational speed-ups—as verified across a broad suite of numerical experiments. 

The paper is organized as follows. In Section~\ref{sec:background}, we introduce the essential concepts of Markov chains and Riemannian optimization required to formulate the problem. Section~\ref{sec:manifold} presents the construction of the specific manifold on which the optimization is performed and details the necessary operations. In Section~\ref{sec:numerical_examples}, we validate the proposed algorithm through numerical experiments on both synthetic test problems and Markov chains arising from the simulation of stochastic differential equations. Finally, Section~\ref{sec:conclusions} summarizes our conclusions and discusses potential directions for future research.

\section{Background and problem formulation}\label{sec:background}

To build the approach we discuss here we need to introduce some concepts regarding discrete-time finite Markov chains, stochastic matrices and Riemannian optimization. We start by recalling the connection between stochastic matrices and Markov chains; for a complete presentation of the topic one can refer to the classic book by Kemeny and Snell\cite{KemenySnell}. We first fix several notations used throughout this work. Given $A,B \in \mathbb{R}^{m\times n}$ we denote by the $A \odot B$ and $ A \oslash B$ the Hadamard entrywise multiplication and division, respectively. Moreover, given a matrix $A\in \mathbb{R}^{m \times n}$, we denote by $A \geq 0$ or $A > 0$ an entrywise nonnegative and an entrywise positive matrix, respectively. Given a vector $\boldsymbol{v}\in \mathbb{R}^n$, we denote by $D_{\boldsymbol{v}} = \operatorname{diag}(\mathbf{v})$ the diagonal matrix, whose diagonal elements are equal to the entries $v_1,\ldots,v_n$ of the vector~$\boldsymbol{v}$, and denote by $\operatorname{diag}(A)$ the main diagonal of its matrix argument.

\begin{definition}
A finite-dimensional discrete-time Markov chain is a discrete-time stochastic process $\{X_\ell\}_{\ell \geq 0}$ defined on a probability space $(\Omega, \mathcal{F}, \mathbb{P})$ that takes values in a finite-dimensional state space $\mathcal{V}$ which, for simplicity, is assumed to be $\mathcal{V}=\{1,\ldots,n\}$. The process satisfies the Markov property, meaning that for every integer $\ell \geq 0$ and any states $i_0, i_1, \dots, i_{\ell+1} \in \mathcal{V}$, the following holds:
\[
\mathbb{P}(X_{\ell+1}=i_{\ell+1} \mid X_0=i_0, \dots, X_\ell=i_\ell) = \mathbb{P}(X_{\ell+1}=i_{\ell+1} \mid X_\ell=i_\ell);
\]
Furthermore, the Markov chain is called homogeneous if
\[
\mathbb{P}(X_{\ell+1}=j \mid X_\ell=i) = p_{ij}, \qquad \forall\, i,j\in\mathcal{V},
\]
where \( p_{ij} \) represents the probability of transitioning from state \( i \) to state \( j \) in one time unit. The matrix $P=(p_{ij})_{ij}$ is called the transition matrix associated with the Markov chain $\{ X_\ell \}_\ell$.
\end{definition}

\begin{definition}
A {row-}stochastic matrix \( P \in \mathbb{R}^{n \times n} \) is a square matrix satisfying the non-negativity condition \( p_{ij} \geq 0 \) for all \( i, j \in \{1,\ldots,n\} \), and the row-sum conditions $P\mathbf{1} = \mathbf{1}$, for $\mathbf{1} = [1,\ldots,1]^\top$.
\end{definition}

{In the sequel, we will use the term stochastic to denote row-stochastic matrices.}

A transition matrix is a stochastic matrix.

\begin{definition}
A vector $\boldsymbol{\pi}^\top = [\pi_1, \pi_2, \dots, \pi_n]$ is called a stationary vector (or stationary distribution) for a Markov chain with transition matrix $P$ if it satisfies two conditions. First, it is a probability vector, meaning that each $\pi_i \geq 0$ and $\boldsymbol{\pi}^\top \mathbf{1}  = 1$. 
Second, it fulfills the invariance relation
\begin{equation}\label{eq:stationary_distribution}
    \boldsymbol{\pi}^\top = \boldsymbol{\pi}^\top P.
\end{equation}
\end{definition}

\begin{definition}
A subset $C$ of the state space $\mathcal{V}$ of a Markov chain is called an ergodic class if it satisfies the following conditions. First, it is closed; that is, if a state $i$ belongs to $C$ and there is a positive probability of transitioning from $i$ to some state $j$, then $j$ must also belong to $C$. Second, $C$ is irreducible, meaning that for any two states $i$ and $j$ in $C$, there exists a positive integer $\ell$ such that the $\ell$-step transition probability $P^\ell(i,j)$ is positive. Finally, the class is aperiodic; that is, for every state $i\in C$, the greatest common divisor of the set $\{ \ell\geq 1: P^\ell(i,i)>0\}$ is one.
\end{definition}

A Markov chain with a single ergodic class admits a unique stationary distribution vector $\boldsymbol{\pi}$. Moreover, if addition  there are no transient states (i.e., no state $i \in \mathcal{V}$ for which the probability of eventually returning to $i$, starting from $i$, is less than one), then $\boldsymbol{\pi}$ is strictly positive.

Conversely, if the Markov chain has multiple ergodic classes along with a (possibly empty) set of transient states, then the initial state influences the long-run average visitation frequencies. Equivalently, there exist different linearly independent stationary distributions satisfying~\eqref{eq:stationary_distribution}. In this general case, the transition matrix $P$ can, up to a suitable permutation, be expressed as the block matrix
\begin{equation}\label{eq:ergodic_decomposition}
    P = \begin{bmatrix}
P_1 & 0   & 0   & \cdots & 0      \\
0   & P_2 & 0   & \cdots & 0      \\
\vdots & \ddots & \ddots & \ddots & \vdots \\
0   & \cdots & 0 & P_E   & 0      \\
P_{T1} & P_{T2} & \cdots & P_{TE} & P_{TT}
\end{bmatrix},
\end{equation}
where the number of stochastic irreducible square blocks $E$ is the number of ergodic classes and $T$ represents the (possibly empty) set of transient states. If $E > 1$, then $P$ admits infinitely many stationary vectors of the form 
\begin{equation}\label{eq:stationary_reducible}
 \boldsymbol{\pi}^\top  = [\alpha_1 \boldsymbol{\pi}_1^\top,\ldots, \alpha_E \boldsymbol{\pi}_E^\top,\mathbf{0}_{|T|}^\top],    
\end{equation}
where $\alpha_i \in [0,1]$, $\sum_{i=1}^E \alpha_i=1$, $\{\boldsymbol{\pi}_i^\top = \boldsymbol{\pi}_i^\top P_i,~\boldsymbol{\pi}_i^\top \mathbf{1}=1\}_{i=1}^{E}$, and $\mathbf{0}_{|T|}$ is the vector of zeros of size the cardinality of the set of transient states. 

Among all the Markov chains, we are interested here in the case of reversible chains.
\begin{definition}
A Markov chain with finite state space \( \mathcal{V} = \{1, \dots, n\} \) and transition matrix \( P \in \mathbb{R}^{n \times n} \) is called reversible if there exists a probability distribution \( \boldsymbol{\pi}^\top = [\pi_1, \dots, \pi_n] \) satisfying \( \pi_i > 0 \) for all \( i \in \mathcal{V} \) and \( \sum_{i \in \mathcal{V}} \pi_i = 1 \), such that the detailed balance condition holds:
\begin{equation}\label{eq:detailed_balance}
    \pi_i P_{ij} = \pi_j P_{ji}, \quad \forall i,j \in \mathcal{V},
\end{equation}
or, in matrix format, 
\[
D_{\boldsymbol{\pi}} P = P^\top D_{\boldsymbol{\pi}}.
\]
\end{definition}

By multiplying on the left the above equation by $\mathbf{1}^\top$ we obtain the equality $\boldsymbol{\pi}^\top P=\boldsymbol{\pi}^\top $, i.e,
$\boldsymbol{\pi}$ is the stationary distribution of $P$.
With an abuse of notation, a stochastic matrix $P$, with stationary vector $\boldsymbol{\pi}$, will be called reversible if it satisfies \eqref{eq:detailed_balance}.
{Note, however, that reversibility does \emph{not} imply $P = P^\top$; if $\boldsymbol{\pi}$ is uniform, then a reversible matrix is also symmetric.}

The reversibility condition~\eqref{eq:detailed_balance} ensures that the probability flow from state \( i \) to state \( j \) is equal to the probability flow from state \( j \) to state \( i \), making the Markov chain behave as if it were in equilibrium. A Markov chain satisfying this condition is time-reversible, meaning that the process appears the same when observed in reverse time.  While this is the usual setting for which reversibility is defined\cite{KemenySnell}%
, for our task we want to consider also the case of a general reducible chain with possibly both multiple ergodic classes and transient states; if transient states are present time-reversibility of the underlying chain is no longer possible, but the detailed balance condition~\eqref{eq:detailed_balance} is still meaningful. In this scenario, the conditions~\eqref{eq:detailed_balance} can be reformulated for each ergodic class individually by considering any vector $\boldsymbol{\pi}$ of the form~\eqref{eq:stationary_reducible} {where $\alpha_i>0$ for $i=1,\ldots,E$}; observe also that~\eqref{eq:detailed_balance} are invariant with respect to any choice of scaling factors $\alpha_i$'s in $\boldsymbol{\pi}$---globally and for the single classes. Then the detailed balance equations~\eqref{eq:detailed_balance} for indices corresponding to the transient states are such that
\[
 \underbrace{\pi_i}_{= 0} P_{ij} = \pi_j \underbrace{P_{ji}}_{= 0}, \quad \forall\, i \in T,\; j \in \mathcal{V}.
\]
Since both $\pi_i$ and $P_{ji}$ are zero for any $i \in T$, these equations are trivially satisfied. This observation indicates that, in the presence of both transient and ergodic classes,
the transient states can be removed from the set of states, and
the reversibility conditions~\eqref{eq:detailed_balance}
can be applied to the remaining states, where $\boldsymbol{\pi}$ is chosen to have all positive entries.

\subsection{Optimization problem}

In this subsection, we provide a formulation of the problem of finding the nearest reversible Markov chain as an optimization problem. Let $A \in \mathbb{R}^{n \times n}$ be a matrix belonging to the set
\[
\mathbb{S}_n^0 = \{ S \in \mathbb{R}^{n \times n} : S\mathbf{1} = \mathbf{1}, \; S \geq 0 \},
\]
and let $\boldsymbol{\pi}$ be its stationary vector, i.e., $\boldsymbol{\pi}^\top A = \boldsymbol{\pi}^\top$, which we may assume without loss of generality being $\boldsymbol{\pi} > 0$, since in cases where transient states are present, one may discard the entries corresponding to these states from $\boldsymbol{\pi}$ in~\eqref{eq:stationary_reducible}, and work on the Markov chain from which these states have been removed, as we will observe in Section~\ref{sec:ergodic_division}.

Given such a matrix $A \in \mathbb{S}_n^0$ with stationary vector $\boldsymbol{\pi}>0$, we address the problem of finding the closest stochastic matrix $X$ in the~set
\begin{equation}\label{eq:feasible-set}
    \mathcal{R} = \Big\{ S \in \mathbb{R}^{n \times n} : S \mathbf{1} = \mathbf{1}, \; \boldsymbol{\pi}^\top S = \boldsymbol{\pi}^\top, \; D_{\boldsymbol{\pi}} S = S^\top D_{\boldsymbol{\pi}}, \; S \geq 0 \Big\} \equiv \Big\{ S \in \mathbb{R}^{n \times n} : S \mathbf{1} = \mathbf{1}, \; D_{\boldsymbol{\pi}} S = S^\top D_{\boldsymbol{\pi}}, \; S \geq 0 \Big\},
\end{equation}
which means solving the optimization problem
\begin{equation}\label{eq:optimization_problem}
    X^\star = \arg \min_{X \in \mathcal{R}}\frac{1}{2} \| X-A \|_F^2,
\end{equation}
where $\| \cdot \|_F$ denotes the Frobenius matrix norm. The set $\mathcal{R}$ collects stochastic reversible matrices, whose stationary distribution vector is equal to $\boldsymbol{\pi}$.

%
%

{The optimization problem~\eqref{eq:optimization_problem} can be reformulated as a standard quadratic programming problem\cite{Nielsen2015483}---namely, one with a quadratic objective function subject to linear equality and inequality constraints---and can therefore be solved by using any suitable quadratic–programming solver\cite{NocedalWright}. Here, instead, we exploit the additional geometric structure of the constraint set in order to develop a smooth optimization framework. Indeed, reversible stochastic matrices are not symmetric in general, and therefore do not form a manifold of symmetric matrices directly. To build a Riemannian formulation amenable to manifold–based numerical methods\cite{MR4533407}, we restrict the feasible set to 
\[
\mathcal{R}_{>} = \mathcal{R} \cap \{ S \in \mathbb{R}^{n\times n} : S >0 \},
\]
and introduce the change of variable
\begin{equation}
\label{eq:change-of-variable}
\hat{S} = D_{\hat{\boldsymbol{\pi}}} S D_{\hat{\boldsymbol{\pi}}}^{-1}, 
\qquad \hat{\boldsymbol{\pi}} = \boldsymbol{\pi}^{\nicefrac{1}{2}}.
\end{equation}
}
The matrix $\hat S$ is symmetric since, from the detailed balance equation \eqref{eq:detailed_balance}, we get
\begin{align*}
    D_{\boldsymbol{\pi}} S D_{\boldsymbol{\pi}}^{-1} &= S^\top,\\[1mm]
    D_{\boldsymbol{\hat \pi}} \, D_{\boldsymbol{\hat \pi}} S D_{\boldsymbol{\hat \pi}}^{-1} \, D_{\boldsymbol{\hat \pi}}^{-1} &= S^\top,\\[1mm]
    \underbrace{D_{\boldsymbol{\hat \pi}} S D_{\boldsymbol{\hat \pi}}^{-1}}_{\hat{S}} &= \underbrace{D_{\boldsymbol{\hat \pi}}^{-1} S^\top D_{\boldsymbol{\hat \pi}}}_{\hat{S}^\top}.
\end{align*}
Moreover, $\boldsymbol{\hat \pi}$ is an eigenvector of $\hat S$ associated with $1$. Indeed, we have
\[
\hat{S}\hat{\boldsymbol{\pi}} = D_{\boldsymbol{\hat \pi}} S D_{\boldsymbol{\hat \pi}}^{-1} \hat{\boldsymbol{\pi}} = D_{\boldsymbol{\hat \pi}} S \mathbf{1} = D_{\boldsymbol{\hat \pi}} \mathbf{1} = {\boldsymbol{\hat \pi}},
\]
which implies that $\hat{\boldsymbol{\pi}}$ is both the right and the left eigenvector of $\hat{S}$ associated with the eigenvalue $1$, with Euclidean norm equal to~1.%

{Therefore, the transformation in~\eqref{eq:change-of-variable} induces a one–to–one correspondence between reversible Markov chains, with positive transition matrix, and positive symmetric matrices with a prescribed Perron eigenvalue and eigenvector. This yields a smooth manifold structure on which Riemannian optimization algorithms can be applied efficiently and with strong theoretical guarantees.}

In particular, the change of variable in \eqref{eq:change-of-variable} transforms the feasible set $\mathcal{R}_{>}$ into the set  
\begin{equation}
\label{eq:manifold}
\mathcal{M}_{\boldsymbol{\pi}} = \Big\{ S \in \mathbb{R}^{n\times n} : S > 0,\; S = S^\top,\; S \hat{\boldsymbol{\pi}} = \hat{\boldsymbol{\pi}} \Big\},  \qquad \hat{\boldsymbol{\pi}} = \boldsymbol{\pi}^{\nicefrac{1}{2}}.
\end{equation}
The functional in \eqref{eq:optimization_problem} can be modified accordingly, employing the same transformation. Indeed, given $X \in \mathcal{R}_{>}$, observe that:
\begin{align*}
    \| X-A\|_F &= \Big\| D_{\hat{\boldsymbol{\pi}}}^{-1} \underbrace{D_{\hat{\boldsymbol{\pi}}} X D_{\hat{\boldsymbol{\pi}}}^{-1}}_{=\hat{X} \in \mathcal{M}_{\boldsymbol{\pi}}} D_{\hat{\boldsymbol{\pi}}} - A \Big\|_F = \left\| D_{\hat{\boldsymbol{\pi}}}^{-1} \hat{X} D_{\hat{\boldsymbol{\pi}}} - A \right\|_F.%
\end{align*}
Thus, we focus on the solution of the following equivalent optimization problem
\begin{equation}\label{eq:equivalent_optimization}
P^{\star} = \arg \min_{\hat{X} \in \mathcal{M}_{\boldsymbol{\pi}}} \frac{1}{2}\|D_{\hat{\boldsymbol{\pi}}}^{-1} \hat{X} D_{\hat{\boldsymbol{\pi}}} - A\|_F^2,
\end{equation}
from which we recover the nearest reversible Markov chain $P$ to $A$ as
\[
P = D_{\boldsymbol{\hat \pi}}^{-1} P^{\star} D_{\boldsymbol{\hat \pi}}.
\]

{
The following convexity properties hold.}

{
\begin{proposition}\label{prop:convexity}
The set $\mathcal{M}_{\boldsymbol{\pi}}$ in~\eqref{eq:manifold} is convex
and the objective functional
\[
f(\hat{X}) = \frac{1}{2}\|D_{\hat{\boldsymbol{\pi}}}^{-1} \hat{X} D_{\hat{\boldsymbol{\pi}}} - A\|_F^2
\]
is strictly convex in the set $\mathcal{M}_{\boldsymbol{\pi}}$.
\end{proposition}
}

\begin{proof}
{
We first observe that  $\mathcal{M}_{\boldsymbol{\pi}}$ is convex. Indeed, the conditions $S = S^\top$ and $S \hat{\boldsymbol{\pi}} = \hat{\boldsymbol{\pi}}$ are linear constraints, hence they define an affine subspace of $\mathbb{R}^{n\times n}$. Moreover, the componentwise positivity constraint $S > 0$ defines a convex set. Therefore, $\mathcal{M}_{\boldsymbol{\pi}}$ is the intersection of an affine subspace with a convex set, and is thus convex.}

{
Let us prove that the functional $f(\hat{X})$ is strictly convex. Since $D_{\hat{\boldsymbol{\pi}}}$ is a diagonal matrix with strictly positive diagonal entries $\hat{\pi}_i > 0$, it is invertible. Hence, the mapping
\[
\mathcal{L} : \hat{X} \mapsto D_{\hat{\boldsymbol{\pi}}}^{-1} \hat{X} D_{\hat{\boldsymbol{\pi}}}
\]
is a linear bijection on $\mathbb{R}^{n\times n}$. In particular, $\mathcal{L}$ is injective.
The functional $f$ can therefore be written as
\[
f(\hat{X}) = \frac{1}{2}\|\mathcal{L}(\hat{X}) - A\|_F^2,
\]
that is, as the squared Frobenius norm of an affine transformation of $\hat{X}$. Since the squared Frobenius norm is strictly convex and $\mathcal{L}$ is injective, it follows that $f$ is strictly convex on $\mathbb{R}^{n\times n}$. 
Restricting a strictly convex function to a convex set preserves strict convexity.
}
\end{proof}

{
As a consequence of Proposition \ref{prop:convexity}, we get that \eqref{eq:equivalent_optimization} is a strictly convex optimization problem over the convex set $\mathcal{M}_{\boldsymbol{\pi}}$.
}

{We observe that the set $\mathcal{M}_{\boldsymbol{\pi}}$ is not closed. In particular, we have that the minimum of $f(\hat{X})$ exists unique in the closed set
 \begin{equation*}%
 \overline{\mathcal{M}_{\boldsymbol{\pi}}} = \Big\{ S \in \mathbb{R}^{n\times n} : S \ge
 0,\; S = S^\top,\; S \hat{\boldsymbol{\pi}} = \hat{\boldsymbol{\pi}} \Big\}, 
 \end{equation*}
 which is the closure of $\mathcal{M}_{\boldsymbol{\pi}}$. We choose to work with the set $\mathcal{M}_{\boldsymbol{\pi}}$, instead of its closure, since in the context of stochastic matrices, the commonly used metric is the Fisher metric\cite{DurastanteMeini,Douik8861409}, which requires the strict positivity. Indeed, there may be examples where the solution of the optimization problem~\eqref{eq:equivalent_optimization} lies on the boundary of the set $\mathcal{M}_{\boldsymbol{\pi}}$. Nevertheless, we provide an example of this framework in Section~\ref{subsec:zero entries}, where the solution $P^{\star}$ of the optimization problem belongs to the boundary of the set $\mathcal{M}_{\boldsymbol{\pi}}$, and our method provides an approximation of the solution in $\mathcal{M}_{\boldsymbol{\pi}}$.
 }

{Throughout our construction, we have assumed access to a \emph{reliable} stationary distribution vector~$\boldsymbol{\pi}$. In practice, however, this vector may itself be corrupted by noise. In such a setting, one may instead observe a perturbed stationary distribution~${\boldsymbol{\pi}} + {\boldsymbol{\delta}}$. Solving the optimization problem~\eqref{eq:equivalent_optimization} for each of these two choices and invoking the equivalence between the formulations yields the corresponding stochastic matrices $P$ and $P + \Delta$, respectively.}

{\begin{proposition}
\label{prop:lower_bound}
    Given $A \in \mathbb{S}^0_n$, $\boldsymbol{\pi}$ and $\boldsymbol{\pi} + \boldsymbol{\delta}$ two probability vectors, consider $P$ and ${P} +\Delta$ the two reversible matrices obtained by solving the optimization problem~\eqref{eq:equivalent_optimization} for the couples $(A,\boldsymbol{\pi})$ and $(A,{\boldsymbol{\pi}} + \boldsymbol{\delta})$. Then, we have 
\[
\begin{split}
\| \Delta \| & \ge \frac{1}{\| (I - P + \mathbf{1}\boldsymbol{\pi}^\top)^{-1} \|} \frac{\|\boldsymbol{\delta}\|}{\| {\boldsymbol{\pi}} +\boldsymbol{\delta} \|}, \\
    \| \Delta^\top \| & \ge \frac{1}{\| (I - P^\top + \boldsymbol{\pi}\mathbf{1}^\top)^{-1} \|} \frac{\|\boldsymbol{\delta}\|}{\| \boldsymbol{\pi} + \boldsymbol{\delta} \|},
\end{split}
\]
for any induced matrix norm.
\end{proposition}}

{\begin{proof}
By detailed balance for \(P\) and \( P + \Delta\), we have
\[
D_{\boldsymbol{\pi}}P = P^\top D_{\boldsymbol{\pi}},\qquad
D_{{\boldsymbol{\pi}}+\boldsymbol{\delta}}(P +\Delta) = (P +\Delta)^\top D_{{\boldsymbol{\pi}}+\boldsymbol{\delta}}.
\]
Subtracting the two equalities gives
\begin{align*}
D_{\boldsymbol{\pi}}P - D_{{\boldsymbol{\pi}}+\boldsymbol{\delta}} (P + \Delta)  & = P^\top D_{\boldsymbol{\pi}} - (P + \Delta)^\top D_{{\boldsymbol{\pi}}+\boldsymbol{\delta}}, \\
D_{{\boldsymbol{\pi}}+\boldsymbol{\delta}} \Delta- \Delta^\top D_{{\boldsymbol{\pi}}+\boldsymbol{\delta}} & = D_{\boldsymbol{\delta}} P - P^\top D_{\boldsymbol{\delta}}, \\
\| {\Delta}^\top \| \geq & \frac{1}{\| (I - P^\top + \boldsymbol{\pi}\mathbf{1}^\top)^{-1} \|} \frac{\|\boldsymbol{\delta}\|}{\| {\boldsymbol{\pi}} + \boldsymbol{\delta}\|},
\end{align*}
having used that $(I - P^\top)^\dagger (I - P^\top) = (I - P^\top + \boldsymbol{\pi}\mathbf{1}^\top)(I - P^\top) = I - \boldsymbol{\pi} \mathbf{1}^\top$; the other inequality is obtained in the same way.
\end{proof}}

\subsection{Optimization on Riemannian manifolds: definitions and tools}
To outline the construction of the optimization algorithm, we begin by recalling the main definitions from Riemannian geometry\cite{Douik8861409,MR4533407,AbsilBook}.
\begin{definition}\label{def:manifold-definition}
Let $\mathcal{E}$ be a linear space of dimension~$d$. A non-empty subset $\mathcal{M}$ of $\mathcal{E}$ is a smooth \emph{embedded manifold} of $\mathcal{E}$ of dimension $n$ if either
\begin{itemize}
    \item $n = d$ and $\mathcal{M}$ is open in $\mathcal{E}$;
    \item $n = d-k$ for some $k \geq 1$ and, for each $x \in \mathcal{M}$, there exists a neighborhood $U$ of $x$ in $\mathcal{E}$ and a smooth function $h : U \rightarrow \mathbb{R}^k$ such that
    \begin{itemize}
        \item If $y \in U$, then $h(y)= 0$ if and only if $y \in \mathcal{M}$; and
        \item $\operatorname{rank} \mathrm{D}h(x) = k$, for $\mathrm{D}h(x)$ the differential of $h$ at $x$.%
    \end{itemize}
\end{itemize}
\end{definition}
The \emph{tangent space} $\mathcal{T}_{x} \mathcal{M}$ at $x \in \mathcal{M}$ can be defined as follows.
\begin{definition}\label{def:tangent_vector}
    Given a manifold $\mathcal{M}$, with embedded space $\mathcal{E}$ and a point $x \in \mathcal{M}$, we define the tangent space at $x$ as
    \[
    \mathcal{T}_x \mathcal{M} = \left\lbrace \gamma'(0) \,| \; \gamma: I \mapsto \mathcal{M} \; \mbox{is smooth curve and } \gamma(0) = x \right\rbrace,
    \]
    where $I$ is an open interval containing $0$. Each element of $\mathcal{T}_x \mathcal{M}$ is called \emph{tangent vector} to $\mathcal{M}$ at the point $x$.

    The \emph{tangent bundle} of the manifold $\mathcal{M}$ is the disjoint union $\mathcal{T} \mathcal{M}$ that assembles all the tangent vectors, i.e., \[\mathcal{T} \mathcal{M} = {\displaystyle \bigsqcup_{x \in \mathcal{M}}} \mathcal{T}_x \mathcal{M}.\]

\end{definition}

%

A Riemannian manifold $\mathcal{M}$ is a smooth manifold endowed with a positive-definite inner product $\langle \xi_x, \eta_x \rangle_x$ on each tangent space $\mathcal{T}_x\mathcal{M}$, for all $\xi_x,\eta_x\in\mathcal{T}_x\mathcal{M}$. This \emph{Riemannian metric} induces the norm
\[
\|\xi_x\|_x = \sqrt{\langle \xi_x, \xi_x\rangle_x}, \quad \forall\, \xi_x \in \mathcal{T}_x\mathcal{M}.
\]
On such a manifold, one can consider the minimization problem
\begin{equation}\label{eq:genopt}
	\min_{x\in\mathcal{M}} f(x),
\end{equation}
where $f\colon\mathcal{M}\to\mathbb{R}$ is a smooth function.

Optimization methods on Riemannian manifolds rely on mapping local information from the tangent spaces $\mathcal{T}_x\mathcal{M}$ back to the manifold $\mathcal{M}$, thereby generating a sequence of iterates that move along curves on the manifold. Different algorithms primarily differ in how they select the new direction in the tangent space. In general, various classes of optimization methods—such as first-order methods, Newton and quasi-Newton methods, and trust-region methods—can be adapted to this setting\cite{AbsilBook}.

To employ these methods, we need the notions of Riemannian gradient and Riemannian Hessian, which we recall here for completeness.

%
%
%
%
%

\begin{definition}\label{def:riemannian-gradient-and-hessian}
	The \emph{Riemannian gradient} of $f$ at $x\in\mathcal{M}$, denoted by $\operatorname{grad}f(x)$, is the unique tangent vector in $\mathcal{T}_x\mathcal{M}$ satisfying
	\[
	\langle \operatorname{grad}f(x), \xi_x \rangle_x = \mathrm{D} f(x)[\xi_x], \quad \forall\, \xi_x \in \mathcal{T}_x\mathcal{M},
	\]
 where $\mathrm{D} f(x)[\xi_x]$ denotes the directional derivative of $f$ at $x$ in the direction $\xi_x$.
 
	The \emph{Riemannian Hessian} of $f$ at $x$, denoted by $\operatorname{hess}f(x)$, is the linear map from $\mathcal{T}_x\mathcal{M}$ to itself defined by
	\[
	\operatorname{hess}f(x)[\xi_x] = \nabla_{\xi_x} \operatorname{grad}f(x), \quad \forall\, \xi_x \in \mathcal{T}_x\mathcal{M},
	\]
	where $\nabla$ denotes the Levi-Civita connection\cite{Lee2018} on $\mathcal{M}$.
\end{definition}

\begin{definition}\label{def:retraction}
	A \emph{retraction} on $\mathcal{M}$ is a smooth map
	\[
	R\colon \mathcal{T}\mathcal{M} = \bigsqcup_{x \in \mathcal{M}} \mathcal{T}_x\mathcal{M} \to \mathcal{M},
	\]
	such that for every $x\in\mathcal{M}$, the restriction $R_x = R|_{\mathcal{T}_x\mathcal{M}}$ satisfies:
	\begin{itemize}
		\item $R_x(0) = x$ (centering),
		\item The curve $\gamma_{\xi_x}(\tau)= R_x(\tau \xi_x)$ has initial velocity $\xi_x$, that is,
		\[
		\left.\frac{d}{d\tau}R_x(\tau \xi_x)\right|_{\tau=0} = \xi_x, \quad \forall\, \xi_x \in \mathcal{T}_x\mathcal{M}.
		\]
	\end{itemize}
\end{definition}

To characterize retractions for embedded manifolds, we recall the following result. %

\begin{theorem}[{Proposition 4.1.2\cite{AbsilBook}}]\label{thm:retraction-on-embedded-manifold}
	Let $\mathcal{M}$ be an embedded manifold of the Euclidean space $\mathcal{E}$ and let $\mathcal{N}$ be another manifold such that
	\[
	\dim(\mathcal{M}) + \dim(\mathcal{N}) = \dim(\mathcal{E}).
	\]
	Assume there exists a diffeomorphism
	\[
	\phi\colon \mathcal{M}\times\mathcal{N}\to\mathcal{E}^*,
	\]
	where $\mathcal{E}^*$ is an open subset of $\mathcal{E}$, and suppose there exists a neutral element $I\in\mathcal{N}$ satisfying
	\[
	\phi(A,I) = A \quad \forall\, A\in\mathcal{M}.
	\]
	Then, the mapping
	\[
	R_x\colon \mathcal{T}_x\mathcal{M}\to\mathcal{M}, \quad \xi_x \mapsto R_x(\xi_x) = \pi_1\Big(\phi^{-1}(x+\xi_x)\Big),
	\]
	where $\pi_1\colon \mathcal{M}\times\mathcal{N}\to\mathcal{M}$ is the projection onto the first component, defines a retraction on $\mathcal{M}$ for all $x\in\mathcal{M}$ and for $\xi_x$ in a neighborhood of $0_x$.
\end{theorem}

In Section~\ref{sec:manifold}, we demonstrate that the set $\mathcal{M}_{\boldsymbol{\pi}}$, defined in~\eqref{eq:manifold} respectively, is an embedded manifold in the sense of Definition~\ref{def:manifold-definition}. We further show that this set can be endowed with an appropriate Riemannian metric, thereby allowing us to reinterpret the optimization problem~\eqref{eq:equivalent_optimization} as one on Riemannian manifolds.

\subsection{Division into ergodic classes}\label{sec:ergodic_division}

{
In the reducible case, we consider a Markov chain that may consist of multiple ergodic classes together with a (possibly empty) set of transient states. We show how the optimization problem in~\eqref{eq:optimization_problem} can be separately addressed on each ergodic class.} 

{As a first stage, we discard the transient states. Transient states are identified by detecting zero entries in $\boldsymbol{\pi}$, where $\boldsymbol{\pi}$ is any vector in the form \eqref{eq:stationary_reducible}, with  
positive coefficients $\alpha_i$, for $i=1,\ldots, E$.
 Such a vector $\boldsymbol{\pi}$ can be computed as the limit of the evolution of an initial probability vector. More precisely, let $\boldsymbol{\pi}^{(0)} > 0$ be a probability vector with strictly positive entries, and consider
\[
\boldsymbol{\pi}^\top = \lim_{k \to \infty} {\boldsymbol{\pi}^{(0)}}^\top A^k .
\]
The strict positivity of the initial distribution ensures that every state is initially assigned nonzero probability. If the chain is reducible, the limiting distribution $\boldsymbol{\pi}$ may depend on the particular choice of $\boldsymbol{\pi}^{(0)}$, since the mass is eventually concentrated on the recurrent classes that are reachable from the initial support. Nevertheless, the set of indices corresponding to the zero entries of $\boldsymbol{\pi}$ uniquely characterizes the transient states, as such states cannot retain positive probability in the limit. From a computational perspective, this stationary vector can be approximated by computing a left eigenvector associated with the eigenvalue $1$ of the transition matrix, for instance by means of the \texttt{eigs} command in \textsc{Matlab}. By selecting an initial vector with strictly positive entries in the iterative eigensolver and normalizing the resulting eigenvector to satisfy the probabilistic constraint, one obtains a stationary distribution whose numerically negligible components (below a prescribed tolerance) identify the transient states in a simple and effective manner.}

{
Once the transient states are identified,
the problem is then restricted to the submatrix $\widetilde{A}$ and the corresponding reduced stationary vector $\widetilde{\boldsymbol{\pi}}$,  obtained by neglecting the transient states.}

{
Then, the reduced chain can be decomposed into ergodic classes via standard graph algorithms such as Tarjan’s method~\cite{MR304178}, and the transition matrix can be written in the form~\eqref{eq:ergodic_decomposition}. Hence, the problem~\eqref{eq:optimization_problem} can be decoupled into independent subproblems, each involving the search for the nearest reversible Markov chain corresponding to the diagonal blocks $\{\widetilde{A}_i\}_{i=1}^{E}$ in the decomposition of $\widetilde{A}$, each of them associated with the corresponding stationary vector $\tilde{\boldsymbol{\pi}}_i$.}

{In this setting, any suitable optimization strategy can be applied separately to each ergodic class in order to compute the closest reversible Markov chain. These computations are independent across classes and can therefore be carried out in an embarrassingly parallel fashion. This approach has the additional advantage of preserving the structure of the original chain, since no artificial couplings between distinct classes are introduced. Moreover, by decomposing the problem into smaller-dimensional subproblems, the overall computational cost is significantly reduced, leading to a substantial speed-up in practical implementations.}

\section{Construction of the manifold and solution of the Riemannian optimization problem}\label{sec:manifold}

We begin by showing that the set $\mathcal{M}_{\boldsymbol{\pi}}$ defined in \eqref{eq:manifold} is an embedded manifold (in the sense of Definition~\ref{def:manifold-definition}) within the linear space $\mathcal{E}\subseteq \mathbb{R}^{n\times n}$. Observe that $\mathcal{M}_{\boldsymbol{\pi}}$ can be seen as an embedded manifold both in $\mathbb{R}^{n\times n}$ and in the set of real symmetric matrices $\mathcal{S}_n= \{ X \in \mathbb{R}^{n \times n} : X = X^\top \}$. Although both choices are admissible, in this work we see the manifold $\mathcal{M}_{\boldsymbol{\pi}}$ as embedded in the set of symmetric matrices $\mathcal{S}_n$, in order to reduce the dimension of the ambient space.

\begin{remark}
    Observe that the set $\mathcal{R}_{>}$ can be seen as an embedded manifold with ambient space $\left\lbrace X \in \mathbb{R}^{n\times n}: D_{\boldsymbol{\pi}}X = X^\top D_{\boldsymbol{\pi}} \right\rbrace$. As a coherence check, we note that both $\mathcal{R}_{>}$ and $\mathcal{M}_{\boldsymbol{\pi}}$ have the same dimension. Specifically, the set $\mathcal{M}_{\boldsymbol{\pi}}$ has dimension $\frac{n(n+1)}{2} - n = \frac{n(n-1)}{2}$, which is the dimension of the space of symmetric matrices minus the $n$ constraints imposed by the eigenvector condition on $\boldsymbol{\hat \pi}$. Similarly, the dimension of $\mathcal{R}_{>}$ can be derived by noting that, for each $S \in \mathcal{R}_{>}$, the condition $\boldsymbol{\pi}^\top S = \boldsymbol{\pi}^\top$ is automatically satisfied once the constraints $D_{\boldsymbol{\pi}} S = S^\top D_{\boldsymbol{\pi}}$ and $S \mathbf{1} = \mathbf{1}$ are enforced; see again~\eqref{eq:feasible-set}. We decide to work in $\mathcal{M}_{\boldsymbol{\pi}}$ instead of in $\mathcal{R}_{>}$ since the former has simpler expressions for the projection and retraction operators.
\end{remark}

We derive an expression for the tangent space, which is essential for formulating the optimization procedure for the minimization problem~\eqref{eq:equivalent_optimization} on the manifold. At a point $S \in \mathcal{M}_{\boldsymbol{\pi}}$, we may identify the tangent space on $\mathcal{M}_{\boldsymbol{\pi}}$ (Definition~\ref{def:tangent_vector}) as follows.
\begin{lemma}
The tangent space to $\mathcal{M}_{\boldsymbol{\pi}}$ at a point $S \in \mathcal{M}_{\boldsymbol{\pi}}$ is 
\begin{equation}
    \mathcal{T}_S \mathcal{M}_{\boldsymbol{\pi}}= \left\lbrace \xi_S \in \mathcal{S}_n : \; \xi_S \boldsymbol{\hat \pi}=0 \right\rbrace . 
\end{equation}
\end{lemma}

\begin{proof}
    Consider a smooth curve $S(t) \in \mathcal{M}_{\boldsymbol{\pi}}$ such that $S(0) = S$, for sufficiently small values of $t$. Then, we have that:
    \begin{align*}
        S(t) \boldsymbol{\hat \pi} = \boldsymbol{\hat \pi}  &\implies \dot{S}(t) \boldsymbol{\hat \pi} = \mathbf{0} \\
          S(t) = S(t)^\top &\implies  \dot{S}(t) = \dot{S}(t)^\top ,
    \end{align*}
which implies that the tangent space satisfies
\begin{equation*}
    \mathcal{T}_S \mathcal{M}_{\boldsymbol{\pi}} \subseteq \left\lbrace \xi_S \in \mathcal{S}_n : \; \xi_S \boldsymbol{\hat \pi}=0 \right\rbrace . 
\end{equation*}
Vice versa, consider $W$ in the set on the right hand side and the smooth curve $\gamma(t) = S + t W$. We observe that the $\gamma(t)$  satisfies $ \gamma(t) = \gamma(t)^\top $ (since $S,W \in \mathcal{S}_n$). Moreover, we have that $\gamma(t) \boldsymbol{\hat \pi} = \boldsymbol{\hat \pi}$ and $\boldsymbol{\hat \pi}^\top \gamma(t) = \boldsymbol{\hat \pi}^\top$. For sufficiently small values of $t$, we also have that $\gamma_{ij}(t) >0$. Since we have that $\gamma(0)=S$ and $\gamma'(0)= W$, we have that $W \in \mathcal{T}_S \mathcal{M}_{\boldsymbol{\pi}}$, by definition of tangent space. Then, we conclude
\begin{equation*}
    \mathcal{T}_S \mathcal{M}_{\boldsymbol{\pi}} \supseteq \left\lbrace \xi_S \in \mathcal{S}_n : \; \xi_S \boldsymbol{\hat \pi}=0 \right\rbrace .\qedhere 
\end{equation*}
\end{proof}

Then, adding on the tangent space $\mathcal{T}_{S}\mathcal{M}_{\boldsymbol{\pi}}$ the Fisher information metric\cite{MR1800071}%
\begin{equation}\label{eq:fisher-metric}
	\begin{split}
		\; \langle \xi_S, \eta_S \rangle_S 
		=  \sum_{i,j=1}^{n} \frac{ (\xi_S)_{ij} (\eta_S)_{ij} }{S_{ij}} = \operatorname{trace}( (\xi_S \oslash S) \eta_S^\top ), \;\forall\, \xi_S,\eta_S \in \mathcal{T}_S \mathcal{M}_{\boldsymbol{\pi}},
	\end{split}
\end{equation}
defines a structure of Riemannian manifold over $\mathcal{M}_{\boldsymbol{\pi}}$.

In order to solve the minimization problem via Riemannian optimization, we need the notion of Riemannian gradient and Riemannian Hessian for a function $f: \mathcal{M}_{\boldsymbol{\pi}} \mapsto \mathbb{R}$. One tool for their derivation is the orthogonal projection $\Pi_{S}$ onto the tangent space
\[
\Pi_{S}\colon \mathcal{S}_n \to \mathcal{T}_S \mathcal{M}_{\boldsymbol{\pi}}.
\]
To compute the orthogonal projection $\Pi_{S}$ onto the tangent space, we first characterize the orthogonal complement of the tangent space.

\begin{lemma}
\label{lem:ort_compl}
    The orthogonal complement of the tangent space $\mathcal{T}_S \mathcal{M}_{\boldsymbol{\pi}}$, with respect to the Fisher metric, has the expression
	\begin{equation}
        \label{eq:orthog_compl}
	    	\mathcal{T}_S^\perp \mathcal{M}_{\boldsymbol{\pi}} = \{ \xi_S^\perp \in \mathcal{S}_n \,:\, \xi_S^\perp = (\boldsymbol{\alpha} \boldsymbol{\hat \pi}^\top + \boldsymbol{\hat \pi} \boldsymbol{\alpha}^\top) \odot S \},
	\end{equation}
	for $\boldsymbol{\alpha} \in \mathbb{R}^{n}$.
\end{lemma}
\begin{proof}
Firstly, we check the inclusion $\supseteq$. Given a $\xi_S \in \mathcal{T}_S \mathcal{M}_{\boldsymbol{\pi}}$ and $Z= (\boldsymbol{\alpha} \boldsymbol{\hat \pi}^\top + \boldsymbol{\hat \pi} \boldsymbol{\alpha}^\top) \odot S \in \mathcal{T}_S^\perp \mathcal{M}_{\boldsymbol{\pi}}$, we get
\begin{align*}
    \left\langle Z, \xi_S \right\rangle_S = \mbox{trace}(\boldsymbol{\alpha} \boldsymbol{\hat \pi}^\top \xi_S^\top + \boldsymbol{\hat \pi} \boldsymbol{\alpha}^\top \xi_S^\top) = \underbrace{\boldsymbol{\hat \pi}^\top \xi_S^\top}_{=0} \boldsymbol{\alpha} + \boldsymbol{\alpha}^\top \underbrace{\xi_S^\top \boldsymbol{\hat \pi}}_{=0} = 0.
\end{align*}
Secondly, we check that both the set at the right hand side and the one at the left hand side in \eqref{eq:orthog_compl} have the same dimension. Indeed, we have that the set at the right-hand-side has dimension $n$, while the orthogonal complement $\mathcal{T}_S ^{\perp} \mathcal{M}_{\boldsymbol{\pi}}$ has dimension equal to $\dim (\mathcal{S}_n) - \dim(\mathcal{T}_S\mathcal{M}_{\boldsymbol{\pi}}) = \frac{n(n+1)}{2} - \left( \frac{n(n+1)}{2} -n \right) = n$.

Observe that since $S$ is symmetric, the elements $\xi_S^{\perp} \in \mathcal{T}_S^{\perp} \mathcal{M}_{\boldsymbol{\pi}}$ are symmetric as well.
\end{proof}

Using Lemma \ref{lem:ort_compl}, we may write the projection from the ambient space $\mathcal{S}_n$ to the tangent space as
\begin{align*}
    \Pi_{S}: \mathcal{S}_n &\mapsto \mathcal{T}_S \mathcal{M}_{\boldsymbol{\pi}} \\
    Z &\mapsto Z - (\boldsymbol{\alpha} \boldsymbol{\hat \pi}^\top + \boldsymbol{\hat \pi} \boldsymbol{\alpha}^\top) \odot S,
\end{align*}
choosing $\boldsymbol{\alpha}$ such that $Z\boldsymbol{\hat \pi} = \left( (\boldsymbol{\alpha} \boldsymbol{\hat \pi}^\top + \boldsymbol{\hat \pi} \boldsymbol{\alpha}^\top) \odot S  \right)\boldsymbol{\hat \pi}$. Such value for the vector $\boldsymbol{\alpha}$ can be derived as
\begin{align*}
    Z\boldsymbol{\hat \pi} = \left( (\boldsymbol{\alpha} \boldsymbol{\hat \pi}^\top + \boldsymbol{\hat \pi} \boldsymbol{\alpha}^\top) \odot S  \right)\boldsymbol{\hat \pi} &= (\boldsymbol{\alpha} \boldsymbol{\hat \pi}^\top \odot S)\boldsymbol{\hat \pi}+ (\boldsymbol{\hat \pi} \boldsymbol{\alpha}^\top \odot S)\boldsymbol{\hat \pi} \\
    &= \operatorname{diag}(S D_{\boldsymbol{\hat \pi}} \boldsymbol{\hat \pi}) \boldsymbol{\alpha} + D_{\boldsymbol{\hat \pi}} S D_{\boldsymbol{\hat \pi}} \boldsymbol{\alpha}.
\end{align*}
The previous formula can be derived combining the following observations. %
Using the properties of the Hadamard product, the term $(\boldsymbol{\alpha} \boldsymbol{\hat \pi}^\top \odot S) \boldsymbol{\hat \pi}$ can be written as
\[
(\boldsymbol{\alpha} \boldsymbol{\hat \pi}^\top\odot S) \boldsymbol{\hat \pi}  = \operatorname{diag} ( \boldsymbol{\alpha} \boldsymbol{\hat \pi}^\top D_{\boldsymbol{\hat \pi}} S^\top ) = \operatorname{diag}(S D_{\boldsymbol{\hat \pi}} \boldsymbol{\hat \pi}) \boldsymbol{\alpha}.
\]
The term $(\boldsymbol{\hat \pi} \boldsymbol{\alpha}^\top \odot S)\boldsymbol{\hat \pi}$ can be rewritten as follows: denote $M= \operatorname{diag}(SD_{\boldsymbol{\hat \pi}} \boldsymbol{\hat \pi})$, for $i=1,\ldots,n$ we have
\begin{align*}
 (Z \boldsymbol{\hat \pi})_{i} = \sum_{j=1}^n Z_{ij} \boldsymbol{\hat \pi}_j = \sum_{j=1}^n M_{ij}\boldsymbol{\alpha}_j + \sum_{j=1}^n \boldsymbol{\hat \pi}_i \boldsymbol{\alpha}_j S_{ij} \boldsymbol{\hat \pi}_j = \sum_{j=1}^n M_{ij}\boldsymbol{\alpha}_j + \sum_{j=1}^n (\boldsymbol{\hat \pi}_i S_{ij} \boldsymbol{\hat \pi}_j )\boldsymbol{\alpha}_j.
\end{align*}

Thus, the vector $\boldsymbol{\alpha}$ can be computed solving the linear system 
\begin{equation}\label{eq:linear_system}
    Z \boldsymbol{\hat \pi} = \left(\operatorname{diag}(S D_{\boldsymbol{\hat \pi}} \boldsymbol{\hat \pi}) + D_{\boldsymbol{\hat \pi}} S D_{\boldsymbol{\hat \pi}} \right) \boldsymbol{\alpha}.
\end{equation}

\begin{remark}\label{rem:matrixA_invertible}
Observe that the linear system~\eqref{eq:linear_system} has a unique solution $\boldsymbol{\alpha}$. To this end, we prove that the matrix 
\begin{equation}\label{eq:matrix_A}
    \mathcal{A} = (\operatorname{diag}(SD_{\boldsymbol{\hat \pi}} \boldsymbol{\hat \pi}) + D_{\boldsymbol{\hat \pi}} S  D_{\boldsymbol{\hat \pi}})
\end{equation}
is invertible. This is equivalent to proving that the scaled matrix 
$$
D_{\boldsymbol{\hat \pi}}^{-1}(\operatorname{diag}(SD_{\boldsymbol{\hat \pi}} \boldsymbol{\hat \pi}) + D_{\boldsymbol{\hat \pi}} S  D_{\boldsymbol{\hat \pi}})D_{\boldsymbol{\hat \pi}} = \operatorname{diag}(SD_{\boldsymbol{\hat \pi}}^2 \mathbf{1}) + S D_{\boldsymbol{\hat \pi}}^2
$$
is invertible. Denote $T= SD_{\boldsymbol{\hat \pi}}^2$. By contradiction, assume that there exists $\mathbf{x} \neq \mathbf{0}$ such that
\[
\operatorname{diag}(T \mathbf{1})\mathbf{x} + T\mathbf{x} =0.
\]
Then, we may write $\operatorname{diag}(T \mathbf{1})^{-1} T \mathbf{x} = -\mathbf{x}$, i.e., $\mathbf{x}$ is an eigenvector corresponding to the eigenvalue $-1$. Since the matrix $P=\operatorname{diag}(T \mathbf{1})^{-1} T$ is positive and stochastic, we have that the eigenvalue $\lambda=1$ is the only eigenvalue whose modulus is equal to the spectral radius\cite{MR544666}, therefore $-1$ cannot be an eigenvalue.
\end{remark}

From a practical standpoint, computing the $\boldsymbol{\alpha}$ needed for the projection involves solving a dense, symmetric linear system whose size matches that of the Markov chain for which we seek the nearest reversible approximation. Given the structure of this system, a direct solver based on Cholesky factorization\cite{GolubAndVanLoan} is well-suited for this task.

The expression for the othogonal projection $\Pi_S$ can be employed in deriving the formulae for the Riemannian gradient and Riemannian Hessian of a smooth function $f: \mathcal{M}_{\boldsymbol{\pi}} \mapsto \mathbb{R}$.

\begin{lemma}\label{lem:riemann-gradient}
  Consider a smooth function $f: \mathcal{M}_{\boldsymbol{\pi}} \mapsto \mathbb{R}$, the Riemannian gradient $\operatorname{grad} f(S)$ has the following expression:
\begin{equation}\label{eq:riemanniangrad-fromeucligrad}
    \operatorname{grad}f(S) = \widetilde{\Pi}_S(\operatorname{Grad}f(S) \odot S),
\end{equation}
where the projection $\widetilde{\Pi}_S: \mathbb{R}^{n\times n} \mapsto \mathcal{T}_S \mathcal{M}_{\boldsymbol{\pi}}$ is defined as $\widetilde{\Pi}_S(Z): = {\Pi}_S\left( \nicefrac{Z+Z^\top}{2}\right)$, and $\operatorname{Grad}\,f(S)$ denotes the euclidean gradient.
\end{lemma}

\begin{proof}
   For the euclidean gradient $\operatorname{Grad} f(S)$ it holds
   \[
   \left\langle \operatorname{Grad} f(S), \xi \right\rangle = \mathrm{D} f(S) [\xi], \quad \forall \xi \in \mathbb{R}^{n\times n},
   \]
where $\left\langle \cdot, \cdot \right\rangle$ denotes the Frobenius inner product on $\mathbb{R}^{n\times n}$. In particular, for each $\xi_S \in \mathcal{T}_S \mathcal{M}_{\boldsymbol{\pi}}$, it holds
    \[
    \left\langle \operatorname{Grad} f(S), \xi_S \right\rangle =   \left\langle \widetilde{\Pi}_S (\operatorname{Grad} f(S) \odot S),  \xi_S \right\rangle_S = \mathrm{D} f(S) [\xi_S],
    \]
    where we employed that the term $(\operatorname{Grad} f(S) \odot S) = \widetilde{\Pi}_S (\operatorname{Grad} f(S) \odot S) + W$, with $W$ belonging to the orthogonal complement of $\mathcal{T}_S \mathcal{M}_{\boldsymbol{\pi}}$ in $\mathbb{R}^{n \times n}$, with respect to the Fisher metric. 
    By definition of Riemannian gradient, we conclude that  $\operatorname{grad}f(S) = \widetilde{\Pi}_S(\operatorname{Grad}f(S) \odot S)$.
\end{proof}

 We can obtain the Riemannian Hessian $\operatorname{hess} f(S)[\xi_S]$ (Definition~\ref{def:riemannian-gradient-and-hessian}) starting from the Euclidean one $\operatorname{Hess} f(S)$ by using the projection operator that we have introduced.
\begin{lemma}\label{lem:riemannian-hessian}
    The Riemannian Hessian $\operatorname{hess} f(S)[\xi_S]$ can be obtained from the Euclidean gradient $\operatorname{Grad} f(S)$ and the Euclidean Hessian $\operatorname{Hess} f(S)$ by using the identity
	\[
	\operatorname{hess} f(S)[\xi_S] = \widetilde{\Pi}_S(\mathrm{D}(\operatorname{grad} f(S))[\xi_S])-\frac{1}{2} \widetilde{\Pi}_S((( \widetilde{\Pi}_S(\operatorname{Grad} f(S))\odot S) \odot \xi_S) \oslash S),
	\]
 where $\mathrm{D}(\operatorname{grad} f(S))[\xi_S]$ can be derived as
 \[
 \mathrm{D}(\operatorname{grad} f(S))[\xi_S] = \dot{\gamma}[\xi_S] - (\dot{\boldsymbol{\alpha}}[\xi_S]  \boldsymbol{\hat \pi}^\top + \boldsymbol{\hat \pi}\dot{\boldsymbol{\alpha}}[\xi_S]^\top)\odot S - (\boldsymbol{\alpha}\boldsymbol{\hat \pi}^\top + \boldsymbol{\hat \pi}\boldsymbol{\alpha}^\top)\odot \xi_S.
 \]
 The quantities $\boldsymbol{\alpha}$ and $ \dot{\boldsymbol{\alpha}}[\xi_S]$, solve respectively the nonsingular linear systems:
 \begin{align*}
     \mathcal{A} \boldsymbol{\alpha}  &= \frac{(\gamma + \gamma^\top)}{2}\boldsymbol{\hat \pi}, \\ %
  \mathcal{A} \dot{\boldsymbol{\alpha}}[\xi_S] & = \mathbf{b} %
  \equiv \frac{(\dot{\gamma}[\xi_S] + \dot{\gamma}[\xi_S]^\top)}{2}\boldsymbol{\hat \pi} - (\mathrm{diag}(\xi_S D_{\boldsymbol{\hat \pi}} \boldsymbol{\hat \pi}) + D_{\boldsymbol{\hat \pi}} \xi_S D_{\boldsymbol{\hat \pi}}) \boldsymbol{\alpha};
 \end{align*}
for $\mathcal{A}$ in~\eqref{eq:matrix_A}, 
and $\dot{\gamma}[\xi_S]$ is computed as
\[
\dot{\gamma}[\xi_S]  = \operatorname{Hess}f(S)[\xi_S]\odot S + \operatorname{Grad}f(S)\odot \xi_S.
\]
\end{lemma}

\begin{proof}
    The expression for the Riemannian Hessian $\operatorname{hess} f (S) [\xi_S]$ depends on the metric---the Fisher metric\cite{DurastanteMeini}~\eqref{eq:fisher-metric} in our case. In particular, we make use of the Koszul formula---Theorem 5.3.1\cite{AbsilBook}---to derive an expression of the Riemannian Hessian via the Levi-Civita connection, as in Definition~\ref{def:riemannian-gradient-and-hessian}. Therefore, in the following, we only highlight the main difference, which consists in the formulae for the derivation of $\mathrm{D}(\operatorname{grad} f(S))[\xi_S]$. Denote $\gamma= \operatorname{Grad} f(S)\odot S$; then, we compute
    \begin{align*}
        \mathrm{D}(\operatorname{grad} f(S))[\xi_S] &=   \mathrm{D}(\widetilde{\Pi}_S(\gamma))[\xi_S] = \mathrm{D}(\gamma - (\boldsymbol{\alpha}\boldsymbol{\hat \pi}^\top + \boldsymbol{\hat \pi}\boldsymbol{\alpha}^\top)\odot S)[\xi_S]  \\
        &= \mathrm{D}(\gamma)[\xi_S] - \mathrm{D}((\boldsymbol{\alpha}\boldsymbol{\hat \pi}^\top + \boldsymbol{\hat \pi}\boldsymbol{\alpha}^\top)\odot S)[\xi_S] \\
        &= \dot{\gamma}[\xi_S] - (\dot{\boldsymbol{\alpha}}[\xi_S]  \boldsymbol{\hat \pi}^\top + \boldsymbol{\hat \pi}\dot{\boldsymbol{\alpha}}[\xi_S]^\top)\odot S - (\boldsymbol{\alpha}\boldsymbol{\hat \pi}^\top + \boldsymbol{\hat \pi}\boldsymbol{\alpha}^\top)\odot \xi_S.
    \end{align*}
The derivative $\dot{\gamma}[\xi_S]$ can be expressed as
\[
\dot{\gamma}[\xi_S] = \mathrm{D}({\operatorname{Grad} f(S)\odot S})[\xi_S]  = \operatorname{Hess}f(S)[\xi_S]\odot S + \operatorname{Grad}\; f(S)\odot \xi_S.
\]
The formula for $\dot{\boldsymbol{\alpha}}[\xi_S]$ can be derived starting from the quantity
\[
\frac{(\gamma + \gamma^\top)}{2}\boldsymbol{\hat \pi} = (\operatorname{diag}(SD_{\boldsymbol{\hat \pi}} \boldsymbol{\hat \pi}) + D_{\boldsymbol{\hat \pi}} S D_{\boldsymbol{\hat \pi}}) \boldsymbol{\alpha},
\]
and then taking the derivative along the direction $\xi_S$, which results in
\begin{align*}
    \frac{(\dot{\gamma}[\xi_S] + \dot{\gamma}[\xi_S]^\top)}{2}\boldsymbol{\hat \pi} = &\; (\operatorname{diag}(\xi_S D_{\boldsymbol{\hat \pi}} \boldsymbol{\hat \pi}) + D_{\boldsymbol{\hat \pi}} \xi_S D_{\boldsymbol{\hat \pi}}) \boldsymbol{\alpha}  +  (\operatorname{diag}(SD_{\boldsymbol{\hat \pi}} \boldsymbol{\hat \pi}) + D_{\boldsymbol{\hat \pi}} S D_{\boldsymbol{\hat \pi}}) \dot{\boldsymbol{\alpha}}[\xi_S].
\end{align*}
Finally, we obtain a formula for $\dot{\boldsymbol{\alpha}}[\xi_S]$ by solving the linear system $\mathcal{A} \dot{\boldsymbol{\alpha}}[\xi_S] = \mathbf{b}$,
where:
\begin{align*}
    \mathcal{A} = &\; (\operatorname{diag}(SD_{\boldsymbol{\hat \pi}} \boldsymbol{\hat \pi}) + D_{\boldsymbol{\hat \pi}} S D_{\boldsymbol{\hat \pi}}), \\ \mathbf{b}  = &\;  \frac{(\dot{\gamma}[\xi_S] + \dot{\gamma}[\xi_S]^\top)}{2}\boldsymbol{\hat \pi} - (\operatorname{diag}(\xi_S D_{\boldsymbol{\hat \pi}} \boldsymbol{\hat \pi}) + D_{\boldsymbol{\hat \pi}} \xi_S D_{\boldsymbol{\hat \pi}}) \boldsymbol{\alpha}, \nonumber
\end{align*}
and $\mathcal{A}$ is non singular by Remark~\ref{rem:matrixA_invertible}.
\end{proof}

In order to conclude the description of the manifold $\mathcal{M}_{\boldsymbol{\pi}}$, we need a characterization for the retraction (Definition~\ref{def:retraction}) from the tangent bundle $\mathcal{T} \mathcal{M}_{\boldsymbol{\pi}}$ (Definition~\ref{def:tangent_vector}) to the manifold $\mathcal{M}_{\boldsymbol{\pi}}$. To this end, we introduce the following lemma, which is a modification of the Sinkhorn-Knopp Theorem\cite{Sinkhorn,MR210731}, or of the DAD-Theorem by Marcus and Newman\cite{MR179182}.

\begin{lemma}
   \label{lemma:DAD_symm}
    Given a symmetric matrix $A\in \mathbb{R}^{n\times n}$, with positive entries, there exists a diagonal matrix $D$ such that
    \begin{equation*}
        DAD \boldsymbol{\hat \pi} = \boldsymbol{\hat \pi}.
    \end{equation*}
\end{lemma}

\begin{proof}
    Consider the symmetric matrix $\hat{A}= D_{\boldsymbol{\hat \pi}} A D_{\boldsymbol{\hat \pi}}$. Then, applying Sinkhorn-Knopp Theorem\cite{Sinkhorn,MR210731}, we get that there exist diagonal matrices $D_1,D_2$ such that
    \[
    D_1 \hat A D_2\mathbf{1} = \boldsymbol{\pi} \; \mbox{and} \; \mathbf{1}^\top D_1 \hat A D_2 = \boldsymbol{\pi}^\top.
    \]
We observe that it also holds
\[
D_2 \hat{A} D_1 \mathbf{1} = \boldsymbol{\pi}, \; \mathbf{1}^\top D_2 \hat{A} D_1 = \boldsymbol{\pi}^\top,
\]
where we used that $\hat{A}=\hat{A}^\top$, since $A$ is symmetric. Since $\hat A_{ij} >0$, we have that the diagonal scales are unique, up to a constant. Then, $D_1= c D_2$. Choosing $D= \sqrt{c}D_1$, we get that $D\hat{A}D \mathbf{1} = \boldsymbol{\pi}$ and $D\hat{A}D$ is symmetric. 

Moreover, we get that
\[
D (D_{\boldsymbol{\hat \pi}} A D_{\boldsymbol{\hat \pi}}) D \mathbf{1} = \boldsymbol{\pi} \implies DAD \boldsymbol{\hat \pi} = \boldsymbol{\hat \pi}.\qedhere
\]
\end{proof}

Using Lemma \ref{lemma:DAD_symm}, we may then define the retraction to the manifold $\mathcal{M}_{\boldsymbol{\pi}}$ as
\begin{align*}
    R: \mathcal{T}\mathcal{M}_{\boldsymbol{\pi}} &\mapsto \mathcal{M}_{\boldsymbol{\pi}},
\end{align*}
with restriction $R_S$ to $\mathcal{T}_S \mathcal{M}_{\boldsymbol{\pi}}$ is given by
\[
R_S( \xi_S) = S + \xi_S,
\]
which is well defined if $\xi_S$ is sufficiently close to $\mathbf{0}_S$ (that is whenever $S > - \xi_S$ entrywise).
Then, by employing Theorem~\ref{thm:retraction-on-embedded-manifold}, we construct the map $\phi$
\begin{align*}
\phi : \mathcal{M}_{\boldsymbol{\pi}} \times \mathbb{R}_{>}^{n} &\mapsto {\overline{\mathcal{S}}_n} \\
\left( S, d \right) &\mapsto \operatorname{diag}(d)\, S\, \operatorname{diag}(d),
\end{align*}
where we define 
\[
{\overline{\mathcal{S}}_n} = \left\lbrace S \in \mathcal{S}_n : \; S_{ij}  > 0 \right\rbrace, \quad \mathbb{R}_{>}^{n} = \left\lbrace x \in \mathbb{R}^{n} \; : x > 0 \right\rbrace.
\]
Note that the sets ${\overline{\mathcal{S}}_n}$ and $\mathbb{R}_{>}^{n}$ are manifold and $\phi$ satisfies the dimensionality property. Indeed, we check that $\dim (\mathcal{M}_{\boldsymbol{\pi}}) = \frac{n(n-1)}{2}$ and 
\[
\dim (\mathcal{M}_{\boldsymbol{\pi}}) + \dim (\mathbb{R}_{>}^{n}) = \frac{n(n-1)}{2} + n = %
\frac{n(n+1)}{2} = \dim ({\overline{\mathcal{S}}_n}).
\]
Finally, we have that for the identity of $\mathbb{R}_{>}^{n}$, we have $\phi(S, \mathbf{1}_{n})= S$. The inverse $\phi^{-1}$ can be constructed via Lemma \ref{lemma:DAD_symm}. 

Therefore, a retraction for $\xi_S$ in the tangent space on $\mathcal{M}_{\boldsymbol{\pi}}$ in sense of Definition~\ref{def:retraction} is given by $\pi_1(\phi^{-1}(S+ \xi_S))$, where $\pi_1:\mathcal{M}_{\boldsymbol{\pi}} \times \mathbb{R}_{>}^{n} \mapsto \mathcal{M}_{\boldsymbol{\pi}}$ is a projection on the first component and the map $\phi^{-1}$ is the identity since
\begin{align*}
    (S+\xi_S) \boldsymbol{\hat \pi} &= S \boldsymbol{\hat \pi} + \xi_S     \boldsymbol{\hat \pi} =  \boldsymbol{\hat \pi},\\
    \boldsymbol{\hat \pi}^\top (S+\xi_S) &= \boldsymbol{\hat \pi}^\top  S + \boldsymbol{\hat \pi}^\top \xi_S = \boldsymbol{\hat \pi}^\top,\\
    (S + \xi_S)  &= (S+\xi_S)^\top.
\end{align*}
The convergence behavior of the Sinkhorn method, as well as the quality of the resulting solutions, depends significantly on the magnitude of the entries in the matrix to which it is applied\cite{MR2399579}. 

One possible precaution to avoid the worsening of the retraction in presence of small modulus entries consists in employing a first order retraction, whose restriction to $\mathcal{T}_S \mathcal{M}_{\boldsymbol{ \pi}}$ is given by
\begin{equation*}
    R_S(\xi_S) = \mathcal{P} \left( S \odot \mathrm{exp}(\xi_S \oslash S) \right), \quad \forall \xi_S \in \mathcal{T}_S \mathcal{M}_{\boldsymbol{ \pi}},
\end{equation*}
where the exponential $\exp(\cdot)$ is performed entrywise and $\mathcal{P}: {\overline{\mathcal{S}}_n} \mapsto \mathcal{M}_{\boldsymbol{ \pi}}$ represents the application of the modified Sinkhorn theorem proposed in Lemma \ref{lemma:DAD_symm}. Several applications of this type of modifications are available in the literature\cite{Douik8861409, DurastanteMeini}.

Nevertheless, in our numerical implementation of the method, to enhance the output quality and ensure that the resulting matrix is symmetric up to machine precision, we combine the previous retraction with the following  proposal.

Even when starting from a non-symmetric matrix $A$, the Sinkhorn variant described in Lemma~\ref{lemma:DAD_symm} yields two diagonal scaling matrices $D_1$ and $D_2$ such that
\[
D_1 \hat{A} D_2 \mathbf{1} = \boldsymbol{\pi}, \quad \mathbf{1}^\top D_1 \hat{A} D_2 = \boldsymbol{\pi}^\top,
\]
where $\hat{A}$ is an appropriately scaled version of $A$. By transposing both equations, we also obtain
\[
D_2 \hat{A}^\top D_1 \mathbf{1} = \boldsymbol{\pi}, \quad \mathbf{1}^\top D_2 \hat{A}^\top D_1 = \boldsymbol{\pi}^\top.
\]
We can then construct the symmetrized matrix
\[
\frac{1}{2} \left( D_1 \hat{A} D_2 + D_2 \hat{A}^\top D_1 \right),
\]
which satisfies $\frac{1}{2} \left( D_1 \hat{A} D_2 + D_2 \hat{A}^\top D_1 \right)\mathbf{1} = \boldsymbol{\pi}$ and is symmetric by construction. In our case, we choose $\hat{A} = D_{\hat{\boldsymbol{\pi}}} A D_{\hat{\boldsymbol{\pi}}}$, so that the final matrix inherits the desired numerical symmetry and maintains the target marginal constraints. 

\begin{remark}
Although not a geometric property of the manifold $\mathcal{M}_{\boldsymbol{\pi}}$, the ability to generate a random point on the manifold is a useful feature for initializing optimization algorithms. For the manifold considered in this work, a random point can be constructed as follows. We begin by generating a matrix $X$ whose entries are given by the absolute value of independent samples from a standard normal distribution with zero mean and unit variance. We then compute its symmetrized form $Y = D_{\hat{\boldsymbol{\pi}}} \frac{(X + X^\top)}{2} D_{\hat{\boldsymbol{\pi}}}$, and apply the generalized Sinkhorn algorithm, as described in Lemma~\ref{lemma:DAD_symm}, to obtain the appropriate diagonal scaling matrices. Applying these scaling matrices to $X$ yields a matrix that lies on the manifold $\mathcal{M}_{\boldsymbol{\pi}}$, thus providing a valid random initialization point for the optimization procedure.
\end{remark}

We have implemented all these operations in the form of a new manifold for the \textsc{Matlab} version of \textsc{Manopt}\cite{manopt}, specifically, given the vector $\boldsymbol{\hat \pi}$ as $\texttt{hp}$, this corresponds to a \lstinline[language=Matlab]!function M = multinomialsymmetricfixedfactory(hp)!, available in the \texttt{GitHub} repository
\href{https://github.com/miryamgnazzo/nearest-reversible}{\texttt{miryamgnazzo\-/nea\-rest\--re\-ver\-si\-ble}}, returning a structure variable \lstinline[language=Matlab]!M! whose fields are the different operations under standard names; e.g., \lstinline[language=Matlab]|M.egrad2rgrad| which implements the conversion of an Euclidean gradient into the Riemannian one in~\eqref{eq:riemanniangrad-fromeucligrad}.

\subsection{Computing the Euclidean gradient and Hessian}

As observed in the construction of the manifold, the Riemannian gradient and Hessian of a smooth function $f:\mathcal{M}_{\boldsymbol{\pi}} \mapsto \mathbb{R}$ can be obtained by a suitable projection of their Euclidean counterparts, %
as described in Lemma~\ref{lem:riemann-gradient} and Lemma~\ref{lem:riemannian-hessian}. %
From a numerical point of view, it may be convenient to employ a closed form for both the Euclidean gradient and Hessian. As anticipated, our numerical implementation of the method relies on the \textsc{Manopt} package, which only requires the computation of the Euclidean Hessian along a direction.

\begin{proposition}[Euclidean gradient and Hessian]
    Let
\[
F(\hat{X}) = \frac{1}{2}\left\|D_{\hat{\boldsymbol{\pi}}}^{-1}\hat{X}D_{\hat{\boldsymbol{\pi}}} - A\right\|_F^2.
\]
Then, the Euclidean gradient of \(F\) is given by
\[
\operatorname{Grad} F(\hat{X}) =D_{\boldsymbol{\pi}}^{-1}\hat X D_{\boldsymbol{\pi}} - D_{\hat{\boldsymbol{\pi}}}^{-1} A D_{\hat{\boldsymbol{\pi}}},
\]
and the Euclidean Hessian along the direction $V \in \mathbb{R}^{n\times n}$ is
\[
\operatorname{Hess} F(\hat{X}) = D_{\boldsymbol{\pi}}^{-1}V D_{\boldsymbol{\pi}},
\]
where we recall that $\boldsymbol{\hat \pi} \odot \boldsymbol{\hat \pi} = \boldsymbol{\pi}$.
\end{proposition}
\begin{proof}
    The Euclidean gradient of $F$ can be derived writing the directional derivative of $F$ at $\hat X$ in the direction $V \in \mathbb{R}^{n \times n}$, that is
    \begin{align*}
        D F(\hat X)[V] &= \lim_{t \to 0} \frac{F(\hat X + tV ) - F(\hat X)}{t} \bigg|_{t=0} \\
        &=\lim_{t \to 0} \frac{1}{t} \left[ F(\hat X) + \frac{t^2}{2}\left\langle D_{\hat{\boldsymbol{\pi}}}^{-1} V D_{\hat{\boldsymbol{\pi}}}, D_{\hat{\boldsymbol{\pi}}}^{-1} V D_{\hat{\boldsymbol{\pi}}} \right\rangle + 
        \left\langle D_{\hat{\boldsymbol{\pi}}}^{-1} t V D_{\hat{\boldsymbol{\pi}}}, D_{\hat{\boldsymbol{\pi}}}^{-1} \hat{X}D_{\hat{\boldsymbol{\pi}}} - A \right\rangle - F(\hat X) \right] \bigg|_{t=0}  \\
        &= \left\langle D_{\hat{\boldsymbol{\pi}}}^{-1} V D_{\hat{\boldsymbol{\pi}}}, D_{\hat{\boldsymbol{\pi}}}^{-1} \hat{X}D_{\hat{\boldsymbol{\pi}}} - A\right\rangle,
    \end{align*}
    from which we get 
\begin{align*}
         D F(\hat X)[V] = \left\langle D_{\hat{\boldsymbol{\pi}}}^{-1} V D_{\hat{\boldsymbol{\pi}}}, D_{\hat{\boldsymbol{\pi}}}^{-1} \hat{X}D_{\hat{\boldsymbol{\pi}}} - A\right\rangle = \mbox{trace}\left( D_{\hat{\boldsymbol{\pi}}} V^\top D_{\hat{\boldsymbol{\pi}}}^{-1}(D_{\hat{\boldsymbol{\pi}}}^{-1} \hat{X}D_{\hat{\boldsymbol{\pi}}} - A)\right) &= \mbox{trace}\left( V^\top D_{\hat{\boldsymbol{\pi}}}^{-1}(D_{\hat{\boldsymbol{\pi}}}^{-1} \hat{X}D_{\hat{\boldsymbol{\pi}}} - A)D_{\hat{\boldsymbol{\pi}}} \right) \\
         &= 
    \left\langle  V, D_{\hat{\boldsymbol{\pi}}}^{-1} (D_{\hat{\boldsymbol{\pi}}}^{-1} \hat{X}D_{\hat{\boldsymbol{\pi}}} - A) D_{\hat{\boldsymbol{\pi}}}\right\rangle.
\end{align*}
Then, we write the Euclidean gradient as
\[
\operatorname{Grad}F(\hat X)= D_{\hat{\boldsymbol{\pi}}}^{-1} (D_{\hat{\boldsymbol{\pi}}}^{-1} \hat{X}D_{\hat{\boldsymbol{\pi}}} - A) D_{\hat{\boldsymbol{\pi}}} = D_{\boldsymbol{\pi}}^{-1}\hat X D_{\boldsymbol{\pi}} - D_{\hat{\boldsymbol{\pi}}}^{-1} A D_{\hat{\boldsymbol{\pi}}}.
\]
Similarly, the Euclidean Hessian can be derived by performing the directional derivative of $\operatorname{Grad}F$. %
The computation of the Hessian at $\hat X$ along the direction $V\in \mathbb{R}^{n \times n}$ is given by:
\begin{align*}
D (\operatorname{Grad} F)(\hat{X})[V] &= \lim_{t \to 0} \frac{\operatorname{Grad} F(\hat X + tV) - \operatorname{Grad} F(\hat X)}{t} \bigg|_{t=0}=D_{\boldsymbol{\pi}}^{-1} V D_{\boldsymbol{\pi}}.\qedhere
\end{align*}
\end{proof}

\subsection{The Riemannian Nearest Reversible Matrix Algorithm}
\label{sec:algorithm}

We summarize here the entire procedure for computing the closest reversible matrix, in the Frobenius norm, to a given stochastic matrix using Riemannian optimization.

Let $A$ be a given stochastic matrix for which we wish to find the closest reversible matrix $P$. We begin by computing its stationary distribution vector~$\boldsymbol{\pi}$, {as in Section~\ref{sec:ergodic_division}. Once the transient states are identified as described in Section~\ref{sec:ergodic_division}, the problem is then restricted to the submatrix $\widetilde{A}$, corresponding to the erogodic classes.}%

At this stage, we may optionally check for the presence of multiple ergodic classes in the reduced Markov chain, which leads to a decomposition of the form~\eqref{eq:ergodic_decomposition}. If such a decomposition is desired, the problem can be decoupled into independent subproblems, {as described in Section~\ref{sec:ergodic_division}.} %

Alternatively, if we assume a single ergodic class or choose not to decompose the problem, we proceed by running the optimization directly on the entire matrix $\widetilde{A}$. With the full second-order geometry of the manifold $\mathcal{M}_{\boldsymbol{\pi}}$ established, as defined in~\eqref{eq:manifold}, we employ a second-order Riemannian optimization method. Specifically, we make use of the \lstinline[language=Matlab]{trustregion} Riemannian algorithm\cite{MR2335248}.
\begin{algorithm}
	\caption{Compute the Nearest Reversible Markov Chain via Riemannian Optimization}
	\label{alg:nearest_reversible}
	\begin{algorithmic}
        \State \textbf{Input:} Stochastic matrix $A$, stationary distribution $\boldsymbol{\pi}$ of $A${, Boolean variable ``recurse\_ergodic''}.
        \State \textbf{Output:} Reversible matrix $P$ closest to $A$ in the Frobenius norm.
		\State Identify transient states by detecting (numerically) zero entries in $\boldsymbol{\pi}$.
		\State Reduce $A$ to obtain the submatrix $\widetilde{A}$ and the corresponding stationary vector $\widetilde{\boldsymbol{\pi}}$.
		\If{{recurse\_ergodic}}
		\State Decompose $\widetilde{A}$ as in~\eqref{eq:ergodic_decomposition} into diagonal blocks $\{\widetilde{A}_i\}_{i=1}^{E}$.
		\State Decompose $\widetilde{\boldsymbol{\pi}}$ into vectors $\{\widetilde{\boldsymbol{\pi}}_i\}_{i=1}^{E}$ corresponding to each block.
		\For{each diagonal block $\widetilde{A}_i$}
		\State Construct the manifold $\mathcal{M}_{\widetilde{\boldsymbol{\pi}}_i}$ as in~\eqref{eq:manifold}.
		\State Solve the optimization problem~\eqref{eq:equivalent_optimization} on $\widetilde{A}_i$ using the Riemannian \texttt{trustregion} method to obtain $X_i$.
		\State Compute the nearest reversible block as 
		\(
		\widetilde{P}_i = D_{\hat{\widetilde{\boldsymbol{\pi}}}_i}^{-1} X_i D_{\hat{\widetilde{\boldsymbol{\pi}}}_i}.
		\)
		\EndFor
		\State Assemble the block-diagonal matrix $\widetilde{P}$ from $\{\widetilde{P}_i\}_{i=1}^{E}$.
		\Else
		\State Construct the manifold $\mathcal{M}_{\widetilde{\boldsymbol{\pi}}}$ as in~\eqref{eq:manifold}.
		\State Solve the optimization problem~\eqref{eq:equivalent_optimization} on $\widetilde{A}$ using the Riemannian \texttt{trustregion} method to obtain $X$.
		\State Compute the nearest reversible matrix as 
		\(
		\widetilde{P} = D_{\hat{\widetilde{\boldsymbol{\pi}}}}^{-1} X D_{\hat{\widetilde{\boldsymbol{\pi}}}}.
		\)
		\EndIf
		\State Reintegrate the transient states into $\widetilde{P}$ to form the reversible matrix $P$.
		\State \Return $P$
	\end{algorithmic}
\end{algorithm}
For clarity and completeness, the steps described above are summarized compactly in Algorithm~\ref{alg:nearest_reversible}. In the accompanying code, this procedure is implemented in the \lstinline[language=Matlab]{function [P, out] = riemannian_nearest_reversible(A, pv, varargin)}.

\begin{example}
To have a qualitative description of the algorithm, we consider the graph\cite{nr} reported in Fig.~\ref{fig:original_graph}, which contains the observed attendance at 14 social events by 18 Southern women\cite{davis1941deep}. To build a Markov chain, we start from its adjacency matrix and normalize its rows by the corresponding row sum. The stochastic matrix  obtained by this procedure is reducible, with one transient class and one ergodic class. As described in Algorithm~\ref{alg:nearest_reversible}, the transient class is identified by the numerically zero entries of the stationary distribution. Our algorithm is applied to the submatrix of corresponding to the ergodic class, then the final reversible matrix is constructed by adding the transient class, according to Algorithm~\ref{alg:nearest_reversible}.

\begin{figure}[htbp]
    \centering
    \subfigure[Original graph with adjacency matrix\label{fig:original_graph}]{\includegraphics[width=0.45\columnwidth]{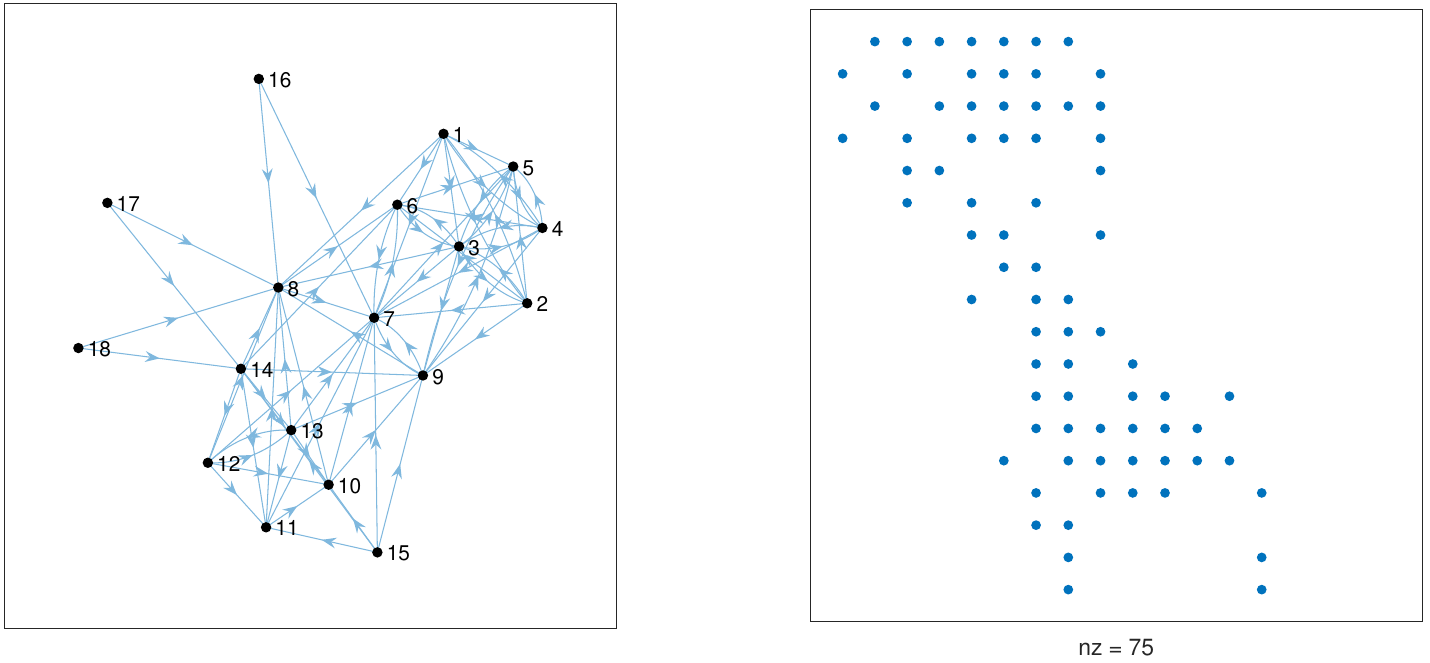}}
    \subfigure[Adjacency graph associated with the nearest reversible matrix obtained with the Riemannian optimization Alg.~\ref{alg:nearest_reversible}]{\includegraphics[width=0.25\columnwidth]{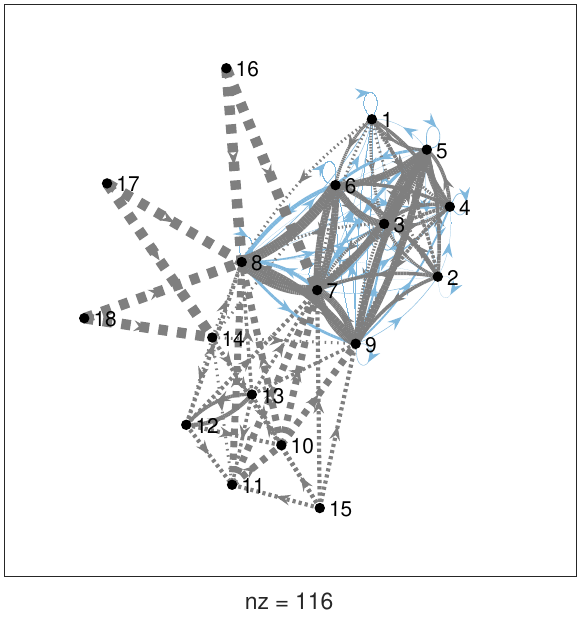}}
    \subfigure[Adjacency graph associated with the nearest reversible matrix obtained with the QP approach\cite{Nielsen2015483}]{\includegraphics[width=0.25\columnwidth]{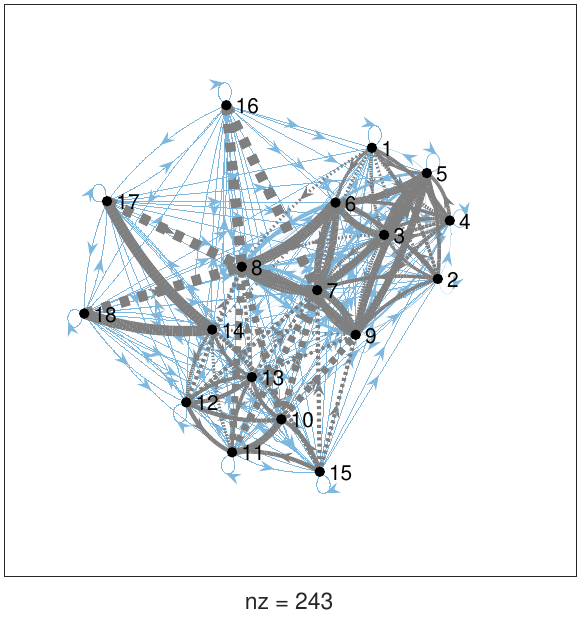}}

    \caption{Example of application of the Algorithm~\ref{alg:nearest_reversible} to the Markov chain built from the graph containing the observed attendance at 14 social events by 18 Southern women\cite{davis1941deep}.}
    \label{fig:small_example}
\end{figure}

The dotted grey lines represent the edges of the original graph, while the solid blue lines represent the new edges. The thickness of the edges is proportional to the corresponding weight.

Observe that our algorithm, by construction, keeps the edges (and their weights) corresponding to the transient states unchanged. Moreover, the new edges are represented by thin lines,  which means that the modification is small in magnitude.
Instead, the algorithm based on the QP approach\cite{Nielsen2015483} adds many new edges, also between nodes belonging to the transient class, by destroying the topology of the original graph. Both algorithm satisfies the detailed balance condition~\eqref{eq:detailed_balance} to machine precision, the QP to \texttt{5.64e-18} and the Riemannian to \texttt{2.9273e-17}. The relative distance in the Frobenius norm between the starting matrix and the one computed by the two algorithms---$\nicefrac{\|A-P\|_F}{\|A\|_F}$---is $\approx 0.312$---with the one obtained via the QP algorithm being larger at the seventh significant digit.
\end{example}

\section{Numerical Examples}\label{sec:numerical_examples}

The proposed algorithm is implemented in \textsc{Matlab}, utilizing the \textsc{Manopt} library\cite{manopt} for Riemannian optimization. The complete source code for the numerical experiments is publicly available in the \texttt{GitHub} repository: 
\href{https://github.com/miryamgnazzo/nearest-reversible}{\texttt{miryamgnazzo\-/nea\-rest\--re\-ver\-si\-ble}}. Notably, the repository includes a \textsc{Manopt}-compatible implementation of the Riemannian manifold~$\mathcal{M}_{\boldsymbol{\pi}}$.

Numerical experiments were conducted on a single node of the Toeplitz cluster, located at the Green Data Center of the University of Pisa. The node is equipped with 256 CPUs (AMD EPYC 7763 64-Core Processor), featuring 2 threads per core, 64 cores per socket, 2 sockets, and a total of 2~TB of RAM.

In Section~\ref{sec:syntethic}, we assess the robustness and effectiveness of the proposed algorithm by applying it to a diverse collection of synthetically generated Markov chains. We also benchmark its performance against the quadratic programming approach\cite{Nielsen2015483}. 

In Section~\ref{sec:application}, we demonstrate the algorithm's practical applicability to Markov chains derived from empirical transition counts of reversible processes. The aim is to recover a reversible transition matrix, which is especially important in settings where approximate transfer operators of the underlying dynamical system result in non-reversible matrices due to numerical errors or insufficient sampling.

\subsection{Synthetic test problems}\label{sec:syntethic}

\begin{figure}[htb]
    \centering
    \subfigure[Uniformly random\label{fig:uniform_random}]{%
        \includegraphics[width=0.23\columnwidth,height=0.23\columnwidth]{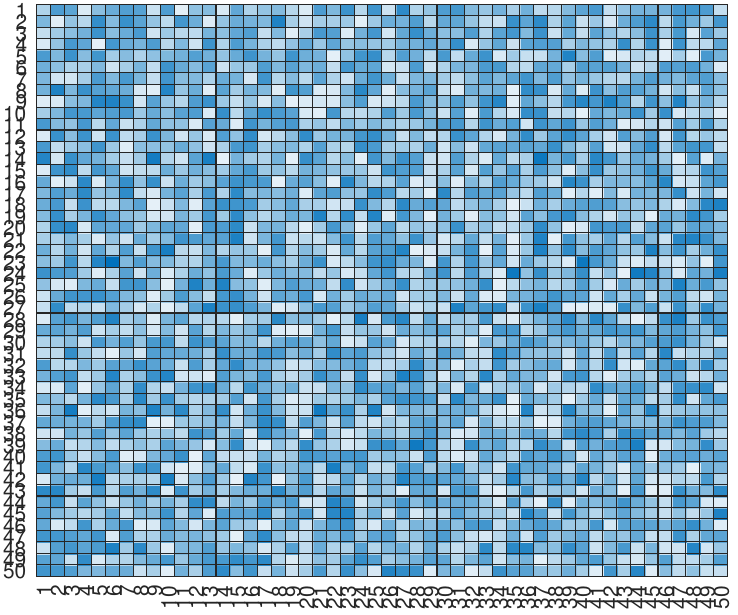}
    }
    \subfigure[Normal random\label{fig:normal_random}]{%
        \includegraphics[width=0.23\columnwidth,height=0.23\columnwidth]{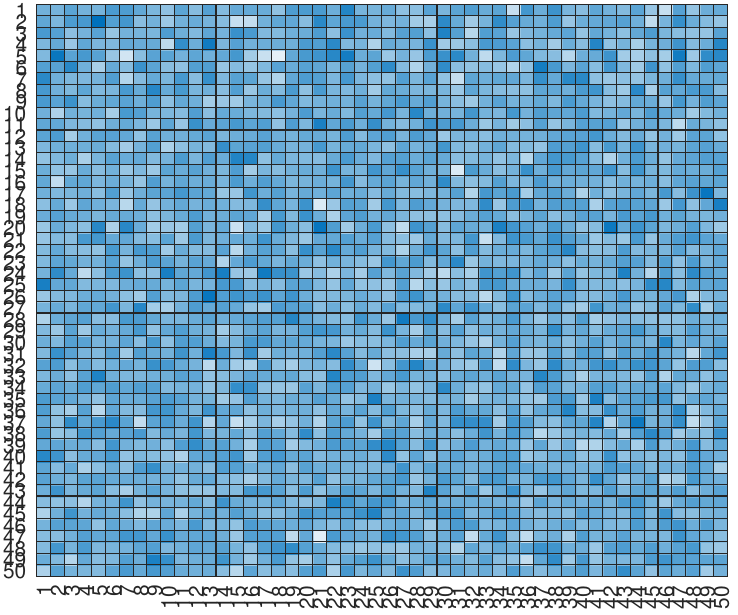}
    }
    \subfigure[Stochastic block model\label{fig:sbm_model}]{%
        \includegraphics[width=0.23\columnwidth,height=0.23\columnwidth]{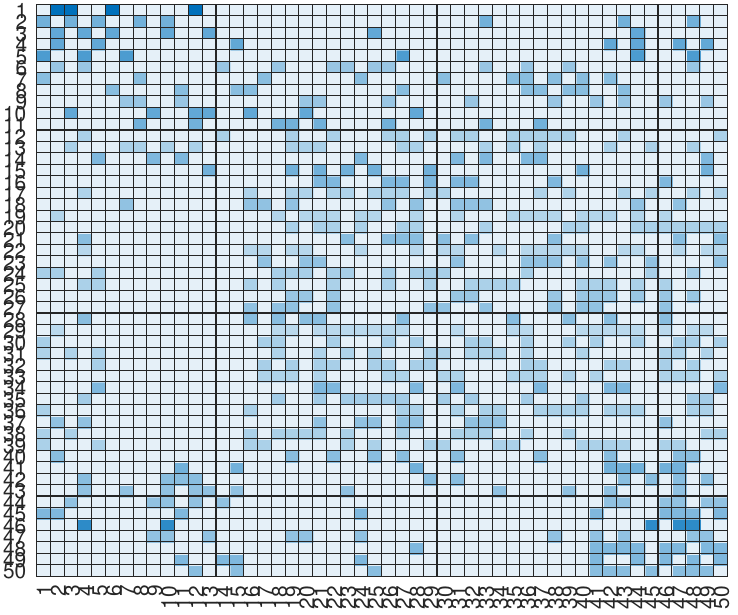}
    }
    \subfigure[Multiple ergodic classes\label{fig:mergodic_model}]{%
        \includegraphics[width=0.23\columnwidth,height=0.23\columnwidth]{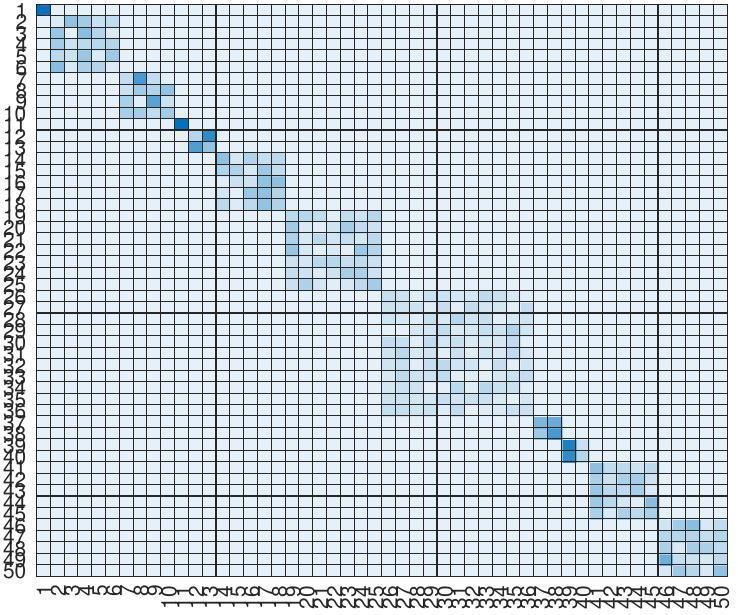}
    }
    
    \caption{Examples of test problems: the stochastic matrices generated from uniformly random entries (see Figure~\ref{fig:uniform_random}), the chain derived from normally distributed random entries (see Figure~\ref{fig:normal_random}), the chain constructed using a stochastic block model (see Figure~\ref{fig:sbm_model}), and the chain with multiple ergodic classes (see Figure~\ref{fig:mergodic_model}).}
    \label{fig:sample_chains}
\end{figure}

To evaluate the proposed algorithm on a range of different problems, we consider the following synthetic test cases:
\begin{description}
    \item[Uniformly random] Generate a matrix \(G\) with entries drawn uniformly from the interval \([0,1]\). Compute the diagonal matrix \(D = \operatorname{diag}(G\mathbf{1})\) and define the stochastic matrix \(A = D^{-1} G\). Figure~\ref{fig:uniform_random} shows an example of a generated matrix using this approach. 
    \item[Normal random] Generate a matrix \(G\) with entries drawn from the normal distribution \(\mathcal{N}(1,1)\). Then, compute \(D = \operatorname{diag}(G\mathbf{1})\) and set \(A = D^{-1} G\). Figure~\ref{fig:normal_random} shows an example of a generated matrix using this approach. 
    \item[Stochastic block model] Generate a matrix \(G\) as the adjacency matrix of a graph obtained from the stochastic block model, and define the associated stochastic matrix \(A\) from the random walk on that graph. This is accomplished using a \textsc{Matlab}-\textsc{Python} interface with the \textsc{NetworkX} library\cite{SciPyProceedings_11} via the command \lstinline[language=Python,basicstyle=\ttfamily,keywordstyle=\rmfamily\bfseries,breaklines]{stochastic_block_model(sizes, p, direct=True)}. To determine the group sizes, we randomly select the number of groups between \(2\) and \(\lfloor n/2 \rfloor\). We then partition \(n\) into \(k\) parts (each at least \(1\)) by choosing \(k-1\) unique cut points from \(\{1,2,\ldots,n-1\}\) so that the differences between consecutive cut points sum to \(n\). The edge density between nodes in group \(r\) and group \(s\) is determined by a randomly generated diagonally dominant stochastic matrix. Figure~\ref{fig:sbm_model} shows an example of a generated matrix using this approach.
    \item[Multiple ergodic classes] We determine class sizes by computing again \(k-1\) unique cut points from \(\{1,2,\ldots,n-1\}\) so that the differences between consecutive cut points sum to \(n\). For each class we generate a uniformly random Markov chain. Figure~\ref{fig:mergodic_model} shows an example of a generated matrix using this approach.
\end{description}

The first two test problems yield generic Markov chains with dense transition matrices, whereas the stochastic block model example represents a nearly decomposable chain that is close to possessing several distinct ergodic classes.

To showcase our results, we employ performance profile graphs\cite{MR1875515}. We compare $3$ solvers: our proposed method, the quadratic programming approach\cite{Nielsen2015483} implemented using \textsc{Matlab}'s \lstinline[language=Matlab]{quadprog} solver, and the same approach using the barrier method from \textsc{Gurobi}\footnote{The solver is commercial software, but it is available at \href{https://www.gurobi.com/}{www.gurobi.com} under a free academic license and has a \textsc{Matlab} interface. We employed \textsc{Gurobi} \texttt{v.11.08}.}. The evaluation is performed over $n_p = 100 + 100 + 100 = 300$ test problems,%
where for each size $n\in \left\lbrace 50,100,200 \right\rbrace$ of the stochastic matrices, we generate $25$ tests drawn from each of the four classes described above.

For each solver $s$ and problem $p$, we consider a set of performance metrics $\{t_{p,s}^{(k)}\}_{k=1}^4$, representing respectively: computational time, relative Frobenius norm distance, $\| D_{\boldsymbol{\pi}} P - P^\top D_{\boldsymbol{\pi}} \|_{\infty}$, and $\| \boldsymbol{\pi}^\top P - \boldsymbol{\pi}^\top \|_\infty$. These metrics quantify both computational efficiency and solution quality.
For each metric $k$, we define the performance ratio:
\[
r_{p,s}^{(k)} = \frac{t_{p,s}^{(k)}}{\min\{ t_{p,s'}^{(k)} : 1 \leq s' \leq 3 \}},
\]
where $r_{p,s}^{(k)} \geq 1$, and a value of 1 indicates that solver $s$ achieved the best performance for problem $p$ under metric $k$.
The performance profile for solver $s$ and metric $k$ is defined as the cumulative distribution function:
\[
\rho_{s}^{(k)}(\tau) = \frac{1}{n_p} \left| \left\{ p : r_{p,s}^{(k)} \leq \tau \right\} \right|,
\]
which represents the fraction of problems for which solver $s$ performs within a factor $\tau$ of the best solver under metric $k$.
In the performance profile graph, $\rho_{s}^{(k)}(\tau)$ is plotted on the vertical axis against $\tau$ on the horizontal axis. The value $\rho_{s}^{(k)}(1)$ reflects the proportion of problems for which solver $s$ was the best. A higher curve indicates better overall performance, while the asymptotic value as $\tau \to \infty$ shows the fraction of problems successfully solved by solver $s$.

In all the following experiments, the test cases are grouped by increasing size and for each size, they are ordered considering one for each class described at the beginning of this section. 

While this grouping does not affect the construction of the performance profiles, it enhances the clarity of the detailed plots. In the following denote the three algorithms as QP-\textsc{Matlab}, Riemannian, and QP-\textsc{Gurobi}. {In all test cases involving Markov chains with multiple ergodic classes, we decompose the problem and apply the chosen optimization procedure independently to each class, as described in Section~\ref{sec:ergodic_division}. For the Riemannian approach introduced in Algorithm~\ref{alg:nearest_reversible}, this class-wise treatment constitutes a natural extension that preserves the structure of the chain and improves computational efficiency. 
For the sake of a fair and meaningful comparison, we adopt the same decomposition strategy for the QP method available in the literature\cite{Nielsen2015483}. Although this strategy was not originally formulated in that context, its use is essential to obtain competitive and computationally reasonable performance across all methods considered.}

{The stopping criteria for QP-\textsc{Matlab} consist of a combined requirement on the relative change of the iterates and on first-order optimality. In particular, the algorithm terminates when the change in the solution falls below $10^{-12}$ and the first-order optimality measure is below $10^{-8}$. For QP-\textsc{Gurobi}, we employ the default parameter settings. Specifically, the tolerance on optimality, \texttt{OptimalityTol}, is set to $10^{-6}$, while the barrier convergence tolerance, \texttt{BarConvTol}, is set to $10^{-8}$. For the Riemannian optimization approach, we use the \texttt{trustregion} algorithm with its default configuration. In this case, convergence is declared when the norm of the Riemannian gradient falls below $10^{-6}$.}

\begin{figure}[htb]
    \centering
    \input{reversibilitylarge}
    \caption{The two plots represent the metric $\| D_{\boldsymbol{\pi}}P - P^\top D_{\boldsymbol{\pi}} \|_\infty$ for $P$, the nearest reversible matrix computed via the different algorithms. The left plot shows the performance profile, while the right plot provides detailed results for each experiment.}
    \label{fig:reversibility_comparison}
\end{figure}

Figure~\ref{fig:reversibility_comparison} compares three algorithms in terms of the reversibility error, computed as \( \|D_{\boldsymbol{\pi}} P - P^\top D_{\boldsymbol{\pi}}\|_\infty \). The left plot shows the performance profile, which illustrates the proportion of problems solved within a factor \(\tau\) of the best-performing method; higher curves indicate better performance across a range of tolerances. The right plot displays the reversibility error on a logarithmic scale for each experiment. The dashed black line represents machine precision (\(2^{-52}\), for \textsc{Matlab} double precision). QP-\textsc{Matlab} and QP-\textsc{Gurobi} consistently achieve errors close to this threshold, while the Riemannian method shows slightly more variability. Nonetheless, all algorithms remain below the machine precision level, indicating that the obtained solutions are numerically reversible.

\begin{figure}[htb]
    \centering
    \input{stationaritylarge}
    \caption{The two plots represent the metric $\| \boldsymbol{\pi}^\top P - \boldsymbol{\pi}^\top \|_\infty$ for $P$, the nearest reversible matrix computed via the different algorithms. The left plot shows the performance profile, while the right plot provides detailed results for each experiment.}
    \label{fig:stationarity_error}
\end{figure}

Figure~\ref{fig:stationarity_error} evaluates how well the computed matrix \( P \) satisfies the stationarity condition, as measured by the metric \( \|{\boldsymbol{\pi}}^\top P - {\boldsymbol{\pi}}^\top\|_\infty \). In the left plot, the performance profile shows all three methods performing well, with QP-\textsc{Matlab} and Riemannian solving nearly all instances accurately. The right plot presents the stationarity error for each experiment on a logarithmic scale. Here, the Riemannian method consistently achieves errors below or near machine precision, outperforming both QP-\textsc{Matlab} and especially QP-\textsc{Gurobi}, the latter of which shows errors reaching up to \(10^{-10}\).

\begin{figure}[htb]
    \centering
    \input{stochasticitylarge}

    \caption{{The two plots represent the metric $\| P \mathbf{1} - \mathbf{1} \|_\infty$ for $P$, the nearest reversible matrix computed via the different algorithms. The left plot shows the performance profile, while the right plot provides detailed results for each experiment.}}
    \label{fig:stochasticity}
\end{figure}

{Figure~\ref{fig:stochasticity} evaluates how accurately the computed transition matrix \( P \) satisfies the stochasticity constraint, quantified by the infinity norm of the residual of the row--sum condition, that is, the deviation of \( P\mathbf{1} \) from \( \mathbf{1} \).  The left panel shows the corresponding performance profile. The Riemannian method clearly dominates, successfully solving almost all instances within very tight tolerances. QP-\textsc{Matlab} attains satisfactory results on a substantial subset of the test set, though with slightly weaker robustness. In contrast, QP-\textsc{Gurobi} exhibits significantly poorer performance, meeting stringent tolerances only in a negligible number of cases. The right panel reports the stochasticity error for each experiment on a logarithmic scale. The Riemannian approach consistently achieves errors at or near machine precision across nearly all instances. QP-\textsc{Matlab} displays a broader spread of errors, yet generally maintains acceptable levels of accuracy. By comparison, QP-\textsc{Gurobi} produces substantially larger violations of the stochasticity constraint, with errors spanning several orders of magnitude and, in extreme cases, becoming exceedingly large. Overall, these results demonstrate the superior numerical stability of the Riemannian formulation in preserving the stochastic structure of the computed Markov chain.}

\begin{figure}[htb]
    \centering
    \input{distancelarge}
    \caption{The two plots represent the relative Frobenius norm distance to the nearest reversible matrix computed via the different algorithms. The left plot shows the performance profile, while the right plot provides detailed results for each experiment.}
    \label{fig:frobenius_distance}
\end{figure}

Figure~\ref{fig:frobenius_distance} illustrates the relative Frobenius norm distance to the nearest reversible matrix, i.e., the quantity $\nicefrac{\| A-P\|_F}{\| A\|_F}$, as computed by each algorithm. As in previous figures, the left plot shows the performance profile. Both QP-\textsc{Matlab} and the Riemannian method solve nearly all problems, retrieving the same distance, while QP-\textsc{Gurobi} performs worse in several cases. The right plot shows the relative distance across experiments. QP-\textsc{Matlab} and Riemannian methods yield consistently low and nearly identical distances. In contrast, QP-\textsc{Gurobi} exhibits significantly higher variability, with several instances showing large deviations, suggesting reduced reliability and failure in multiple cases.

\begin{figure}[htb]
    \centering
    \input{timelarge}
    \caption{The two plots represent the overall computation time (in seconds) to obtain the nearest reversible matrix using the different algorithms. The left plot shows the performance profile, while the right plot provides detailed results for each experiment.}
    \label{fig:time_taken}
\end{figure}

{The final metric considered is computational efficiency, illustrated in Figure~\ref{fig:time_taken}. The left plot presents the performance profile based on runtime. The Riemannian method is the fastest on the majority of instances, followed by QP-\textsc{Matlab}. QP-\textsc{Gurobi} is generally the slowest and exhibits significant variability in execution time. The right plot, shown on a logarithmic scale, reports the runtime per experiment. The Riemannian algorithm completes in under one second for a substantial fraction (approximately one half, and at most about two thirds) of the test cases, whereas the runtime of QP-\textsc{Gurobi} increases markedly with problem size. The growing cost of solving quadratic programs limits the scalability of the QP-based approaches beyond the tested range of stochastic matrix sizes.}

{Overall, the comparison demonstrates that the Riemannian algorithm offers the best balance between accuracy and efficiency for computing the nearest reversible matrix. In terms of reversibility, stationarity and being stochastic (Figures~\ref{fig:reversibility_comparison}, \ref{fig:stationarity_error} and~\ref{fig:stochasticity}), all three methods attain errors at or below machine precision, with the Riemannian approach exhibiting the smallest variability. For the Frobenius-norm distance (Figure~\ref{fig:frobenius_distance}), both the Riemannian and QP-\textsc{Matlab} methods achieve consistently minimal relative distances, whereas QP-\textsc{Gurobi} shows larger deviations and occasional failures. %
Finally, the runtime profiles reported in Figure~\ref{fig:time_taken} clearly demonstrate the efficiency and scalability of the Riemannian approach. The method solves approximately $80\%$ of the test instances within a factor close to the best observed runtime, and more than $95\%$ within a modest performance ratio. In absolute terms, a large fraction of the problems are solved within a few seconds, and nearly all instances are completed within a few minutes, with only a single extreme outlier. By contrast, the QP-based solvers exhibit a markedly slower growth in their performance profiles. Although QP-\textsc{Matlab} performs competitively on a limited subset of small-scale problems, its computational cost increases rapidly with the problem dimension, requiring tens or even hundreds of seconds for a significant portion of the test set. The situation is even more pronounced for QP-\textsc{Gurobi}, whose runtime deteriorates substantially on larger instances, with several experiments requiring hundreds of seconds and the worst case extending to several thousand seconds. The few instances in which QP-\textsc{Matlab} appears faster typically correspond to problems where the decomposition into ergodic classes yields subproblems of very small dimension; in such cases, the fixed overhead associated with the Riemannian optimization framework becomes dominant relative to the actual solution time. Overall, these results confirm that the Riemannian algorithm provides a highly favorable balance between accuracy and computational effort. It consistently achieves near-optimal performance ratios while maintaining substantially lower runtimes, thereby emerging as the method of choice for moderate- and medium-scale problems.}

\subsection{Direct estimation of Markov chains from transition counts}\label{sec:application}

Consider a physical system whose continuous dynamics are modeled by the Langevin stochastic differential equation (SDE)\cite{PhysRevE.93.033307}
\begin{equation}\label{eq:sde}
    \dot{x} = -\frac{\partial U(x)}{\partial x} + \xi(t), \qquad x \in \mathbb{R},
\end{equation}
where $\xi(t)$ is a Gaussian white noise process satisfying
\[
\langle \xi(t)\,\xi(s) \rangle = \sigma^2\,\delta(t-s),
\]
with $\langle \cdot \rangle$ denoting the ensemble average, $\delta(\cdot)$ the Dirac delta function and $U(x)$ a potential function. To simulate the dynamics, one typically employs the Euler--Maruyama scheme\cite{Kloeden1992} to discretize the time evolution using a finite time step~$\Delta t$, on a certain time interval $\left[T_0,T_N\right]$. This yields the update rule
\begin{equation}\label{eq:Euler--Maruyama}
    x(t+\Delta t) = x(t) - \frac{\partial U(x(t))}{\partial x}\,\Delta t + \sigma\,\sqrt{\Delta t}\,\eta(t),
\end{equation}
where $\eta(t)$ represents a sequence of independent standard Gaussian random variables (i.e., with zero mean and unit variance). The discretized trajectory $\{X_{t_n}\}_{n=0}^{N}$, where $T_0 = t_0< t_1< \ldots < t_N=T_N$, thereby serves as an approximation to the underlying continuous-time process.

Next, the continuous state space is partitioned into a collection of $M$ disjoint regions $\{\Omega_i\}_{i=1}^{M}$, which define the discrete states for a Markov model; see an example in the horizontal dashed lines in Figure~\ref{fig:euler-mauryama-example}. The \emph{count matrix} $C$ is constructed by recording the number of transitions between these regions; specifically, the entry $C_{ij}$ counts the number of observed transitions from $\Omega_i$ to $\Omega_j$ within the discretized trajectory. An estimate of the transition probability matrix $A$ is then obtained through row normalization:
\[
A_{ij} = \frac{C_{ij}}{\sum_{k=1}^{M} C_{ik}},
\]
i.e., for every state $\Omega_i \in \mathcal{V}$, the condition $\sum_{j=1}^{M} A_{ij} = 1$ is satisfied.

Although the continuous dynamics are assumed to be reversible, the finite sampling inherent in both the Euler--Maruyama integration and the construction of $C$ may result in an estimated transition probability matrix that does not strictly fulfill the detailed balance condition~\eqref{eq:detailed_balance}. Such deviations necessitate further processing—typically via some kind of projection/reversibilization\cite{Choi_Wolfer} or constrained optimization\cite{Nielsen2015483}—to recover a reversible Markov chain model.

\begin{figure}[htbp]
    \centering
    \input{countestimate}
    \caption{Euler--Maruyama simulation~\eqref{eq:Euler--Maruyama} of an SDE of type~\eqref{eq:sde} where both the potential $U(x) = a + b \cos(x) + c \cos^2(x) + d \cos^3(x)$ and the trajectory $x(t)$ are $2\pi$-periodic, hence the state space is $[0,2\pi]$ and is partitioned in $M = 30$ bins $\{\Omega_i\}_{i=1}^{30}$.}
    \label{fig:euler-mauryama-example}
\end{figure}

\subsubsection{Butane potential} Consider the potential function \( U(x) = a + b \cos(x) + c \cos^2(x) + d \cos^3(x) \), which approximates the butane potential energy profile\cite{Nielsen2015483}. In this expression, \( x \) denotes the central dihedral angle, and the parameters are given by \( a = 2.0567 \), \( b = -4.0567 \), \( c = 0.3133 \), and \( d = 6.4267 \).

Since \( U(x) \) is periodic with period \( 2\pi \), we consider the interval \([0, 2\pi]\), partitioned as illustrated in Figure~\ref{fig:euler-mauryama-example} with $M=30$. Discrete trajectories are generated using the Euler--Maruyama scheme~\eqref{eq:Euler--Maruyama}, with step size \( \Delta t = 10^{-3} \), noise intensity \( \sigma = 1 \), and total number of steps \( n = 5 \times 10^8 \).

We compare two algorithms: the quadratic programming-based method using the \lstinline[language=Matlab]{quadprog} solver\cite{Nielsen2015483} and the Riemannian solver described in Algorithm~\ref{alg:nearest_reversible} with the option to detect the possible presence of different ergodic classes. Because the transition matrix depends on the specific simulation of the stochastic differential equation~\eqref{eq:sde} via the Euler--Maruyama method~\eqref{eq:Euler--Maruyama}, the experiment is repeated \( s = 50 \) times. Results are averaged using the following evaluation metric:%
solution time, satisfaction of the detailed balance equations~\eqref{eq:detailed_balance} in the infinity norm, relative distance from reversibility as estimated by both algorithms, how well \( \boldsymbol{\pi} \) represents the stationary distribution of the optimized matrix, and adherence to the stochasticity condition.

\begin{figure}[htb]
    \centering

    \subfigure[Execution time (s)]{\begin{tikzpicture}

\begin{axis}[%
width=0.238\columnwidth,
height=0.188\columnwidth,
at={(0\columnwidth,0\columnwidth)},
scale only axis,
xmin=0.5,
xmax=2.5,
xtick={1,2},
xticklabels={{QP-\textsc{Matlab}},{Riemannian}},
ymin=0.2837215,
ymax=1.5,
axis background/.style={fill=white},
]
\addplot [color=black, dashed, forget plot]
  table[row sep=crcr]{%
1	0.354427\\
1	0.361347\\
};
\addplot [color=black, dashed, forget plot]
  table[row sep=crcr]{%
2	0.853577\\
2	1.082222\\
};
\addplot [color=black, dashed, forget plot]
  table[row sep=crcr]{%
1	0.34031\\
1	0.345782\\
};
\addplot [color=black, dashed, forget plot]
  table[row sep=crcr]{%
2	0.542719\\
2	0.649454\\
};
\addplot [color=black, forget plot]
  table[row sep=crcr]{%
0.925	0.361347\\
1.075	0.361347\\
};
\addplot [color=black, forget plot]
  table[row sep=crcr]{%
1.925	1.082222\\
2.075	1.082222\\
};
\addplot [color=black, forget plot]
  table[row sep=crcr]{%
0.925	0.34031\\
1.075	0.34031\\
};
\addplot [color=black, forget plot]
  table[row sep=crcr]{%
1.925	0.542719\\
2.075	0.542719\\
};
\addplot [color=blue, forget plot]
  table[row sep=crcr]{%
0.85	0.345782\\
0.85	0.354427\\
1.15	0.354427\\
1.15	0.345782\\
0.85	0.345782\\
};
\addplot [color=blue, forget plot]
  table[row sep=crcr]{%
1.85	0.649454\\
1.85	0.853577\\
2.15	0.853577\\
2.15	0.649454\\
1.85	0.649454\\
};
\addplot [color=red, forget plot]
  table[row sep=crcr]{%
0.85	0.349262\\
1.15	0.349262\\
};
\addplot [color=red, forget plot]
  table[row sep=crcr]{%
1.85	0.733241\\
2.15	0.733241\\
};
\addplot [color=black, only marks, mark=+, mark options={solid, draw=red}, forget plot]
  table[row sep=crcr]{%
1	0.373562\\
1	0.382581\\
1	0.387541\\
1	0.399524\\
1	0.41492\\
1	0.418484\\
1	0.447431\\
1	0.529054\\
1	0.662259\\
};
\addplot [color=black, only marks, mark=+, mark options={solid, draw=red}, forget plot]
  table[row sep=crcr]{%
2	1.196567\\
2	1.27711\\
2	1.296496\\
2	1.835295\\
};
\end{axis}
\end{tikzpicture}%
}
    \subfigure[Relative distance $\nicefrac{\|A-P\|_F}{\|A\|_F}$]{\begin{tikzpicture}

\begin{axis}[%
width=0.238\columnwidth,
height=0.188\columnwidth,
at={(0\columnwidth,0\columnwidth)},
scale only axis,
xmin=0.5,
xmax=2.5,
xtick={1,2},
xticklabels={{QP-\textsc{Matlab}},{Riemannian}},
ymode=log,
ymin=1e-06,
ymax=0.01,
yminorticks=true,
axis background/.style={fill=white},
]
\addplot [color=black, dashed, forget plot]
  table[row sep=crcr]{%
1	0.00441814798074695\\
1	0.00827598981220118\\
};
\addplot [color=black, dashed, forget plot]
  table[row sep=crcr]{%
2	0.00441742940196577\\
2	0.00827540865513824\\
};
\addplot [color=black, dashed, forget plot]
  table[row sep=crcr]{%
1	0.000104856430525662\\
1	0.00167622693742801\\
};
\addplot [color=black, dashed, forget plot]
  table[row sep=crcr]{%
2	1.1028615257348e-06\\
2	0.00167227815030066\\
};
\addplot [color=black, forget plot]
  table[row sep=crcr]{%
0.925	0.00827598981220118\\
1.075	0.00827598981220118\\
};
\addplot [color=black, forget plot]
  table[row sep=crcr]{%
1.925	0.00827540865513824\\
2.075	0.00827540865513824\\
};
\addplot [color=black, forget plot]
  table[row sep=crcr]{%
0.925	0.000104856430525662\\
1.075	0.000104856430525662\\
};
\addplot [color=black, forget plot]
  table[row sep=crcr]{%
1.925	1.1028615257348e-06\\
2.075	1.1028615257348e-06\\
};
\addplot [color=blue, forget plot]
  table[row sep=crcr]{%
0.85	0.00167622693742801\\
0.85	0.00441814798074695\\
1.15	0.00441814798074695\\
1.15	0.00167622693742801\\
0.85	0.00167622693742801\\
};
\addplot [color=blue, forget plot]
  table[row sep=crcr]{%
1.85	0.00167227815030066\\
1.85	0.00441742940196577\\
2.15	0.00441742940196577\\
2.15	0.00167227815030066\\
1.85	0.00167227815030066\\
};
\addplot [color=red, forget plot]
  table[row sep=crcr]{%
0.85	0.00238779248211701\\
1.15	0.00238779248211701\\
};
\addplot [color=red, forget plot]
  table[row sep=crcr]{%
1.85	0.00238333088372491\\
2.15	0.00238333088372491\\
};
\addplot [color=black, only marks, mark=+, mark options={solid, draw=red}, forget plot]
  table[row sep=crcr]{%
1	0.00882352887368383\\
};
\addplot [color=black, only marks, mark=+, mark options={solid, draw=red}, forget plot]
  table[row sep=crcr]{%
2	0.00882265615327332\\
};
\end{axis}
\end{tikzpicture}%
}
    \subfigure[\label{fig:butane:detailed}Detailed balance equation $\|D_{\boldsymbol{\pi}} P - {P}^\top D_{\boldsymbol{\pi}}\|_\infty$]{\begin{tikzpicture}

\begin{axis}[%
width=0.238\columnwidth,
height=0.188\columnwidth,
at={(0\columnwidth,0\columnwidth)},
scale only axis,
unbounded coords=jump,
xmin=0.5,
xmax=2.5,
xtick={1,2},
xticklabels={{QP-\textsc{Matlab}},{Riemannian}},
ymode=log,
ymin=1e-19,
ymax=1e-15,
yminorticks=true,
axis background/.style={fill=white},
]
\addplot [color=black, dashed, forget plot]
  table[row sep=crcr]{%
1	3.47004381446436e-18\\
1	5.20511258685007e-18\\
};
\addplot [color=black, dashed, forget plot]
  table[row sep=crcr]{%
2	6.93911476114401e-18\\
2	1.04086517885295e-17\\
};
\addplot [color=black, dashed, forget plot]
  table[row sep=crcr]{%
1	2.17186661282243e-19\\
1	5.42191455724674e-19\\
};
\addplot [color=black, dashed, forget plot]
  table[row sep=crcr]{%
2	1.73483256042009e-18\\
2	3.46994946417811e-18\\
};
\addplot [color=black, forget plot]
  table[row sep=crcr]{%
0.925	5.20511258685007e-18\\
1.075	5.20511258685007e-18\\
};
\addplot [color=black, forget plot]
  table[row sep=crcr]{%
1.925	1.04086517885295e-17\\
2.075	1.04086517885295e-17\\
};
\addplot [color=black, forget plot]
  table[row sep=crcr]{%
0.925	2.17186661282243e-19\\
1.075	2.17186661282243e-19\\
};
\addplot [color=black, forget plot]
  table[row sep=crcr]{%
1.925	1.73483256042009e-18\\
2.075	1.73483256042009e-18\\
};
\addplot [color=blue, forget plot]
  table[row sep=crcr]{%
0.85	5.42191455724674e-19\\
0.85	3.47004381446436e-18\\
1.15	3.47004381446436e-18\\
1.15	5.42191455724674e-19\\
0.85	5.42191455724674e-19\\
};
\addplot [color=blue, forget plot]
  table[row sep=crcr]{%
1.85	3.46994946417811e-18\\
1.85	6.93911476114401e-18\\
2.15	6.93911476114401e-18\\
2.15	3.46994946417811e-18\\
1.85	3.46994946417811e-18\\
};
\addplot [color=red, forget plot]
  table[row sep=crcr]{%
0.85	1.73512239675191e-18\\
1.15	1.73512239675191e-18\\
};
\addplot [color=red, forget plot]
  table[row sep=crcr]{%
1.85	5.20431560424951e-18\\
2.15	5.20431560424951e-18\\
};
\addplot [color=black, only marks, mark=+, mark options={solid, draw=red}, forget plot]
  table[row sep=crcr]{%
2	1.21431687569308e-17\\
2	1.21433928276298e-17\\
2	1.3877847784873e-17\\
};
\addplot [color=black, dashed, forget plot]
  table[row sep=crcr]{%
0.5	2.22044604925031e-16\\
2.5	2.22044604925031e-16\\
};
\end{axis}
\end{tikzpicture}%
}

    \subfigure[\label{fig:butane:stationary}Stationary vector $\|\boldsymbol{\pi}^\top P - \boldsymbol{\pi}^\top\|_\infty$]{\begin{tikzpicture}

\begin{axis}[%
width=0.238\columnwidth,
height=0.188\columnwidth,
at={(0\columnwidth,0\columnwidth)},
scale only axis,
xmin=0.5,
xmax=2.5,
xtick={1,2},
xticklabels={{QP-\textsc{Matlab}},{Riemannian}},
ymode=log,
ymin=1e-16,
ymax=1e-11,
yminorticks=true,
axis background/.style={fill=white},
]
\addplot [color=black, dashed, forget plot]
  table[row sep=crcr]{%
1	9.52793399733309e-13\\
1	1.77174941384806e-12\\
};
\addplot [color=black, dashed, forget plot]
  table[row sep=crcr]{%
2	3.21964677141295e-15\\
2	6.7168492989822e-15\\
};
\addplot [color=black, dashed, forget plot]
  table[row sep=crcr]{%
1	9.6936347837584e-14\\
1	3.86746190628173e-13\\
};
\addplot [color=black, dashed, forget plot]
  table[row sep=crcr]{%
2	1.11022302462516e-16\\
2	5.55111512312578e-16\\
};
\addplot [color=black, forget plot]
  table[row sep=crcr]{%
0.925	1.77174941384806e-12\\
1.075	1.77174941384806e-12\\
};
\addplot [color=black, forget plot]
  table[row sep=crcr]{%
1.925	6.7168492989822e-15\\
2.075	6.7168492989822e-15\\
};
\addplot [color=black, forget plot]
  table[row sep=crcr]{%
0.925	9.6936347837584e-14\\
1.075	9.6936347837584e-14\\
};
\addplot [color=black, forget plot]
  table[row sep=crcr]{%
1.925	1.11022302462516e-16\\
2.075	1.11022302462516e-16\\
};
\addplot [color=blue, forget plot]
  table[row sep=crcr]{%
0.85	3.86746190628173e-13\\
0.85	9.52793399733309e-13\\
1.15	9.52793399733309e-13\\
1.15	3.86746190628173e-13\\
0.85	3.86746190628173e-13\\
};
\addplot [color=blue, forget plot]
  table[row sep=crcr]{%
1.85	5.55111512312578e-16\\
1.85	3.21964677141295e-15\\
2.15	3.21964677141295e-15\\
2.15	5.55111512312578e-16\\
1.85	5.55111512312578e-16\\
};
\addplot [color=red, forget plot]
  table[row sep=crcr]{%
0.85	6.66799948589869e-13\\
1.15	6.66799948589869e-13\\
};
\addplot [color=red, forget plot]
  table[row sep=crcr]{%
1.85	1.27675647831893e-15\\
2.15	1.27675647831893e-15\\
};
\addplot [color=black, only marks, mark=+, mark options={solid, draw=red}, forget plot]
  table[row sep=crcr]{%
1	1.84868786945458e-12\\
1	2.08655315248052e-12\\
1	3.55976359500687e-12\\
1	3.78758135965995e-12\\
};
\addplot [color=black, only marks, mark=+, mark options={solid, draw=red}, forget plot]
  table[row sep=crcr]{%
2	7.52176099183544e-15\\
2	7.71605002114484e-15\\
2	7.93809462606987e-15\\
2	8.27116153345742e-15\\
2	3.37230243729891e-14\\
2	2.5862645358643e-13\\
2	6.20392626160537e-13\\
};
\addplot [color=black, dashed, forget plot]
  table[row sep=crcr]{%
0.5	2.22044604925031e-16\\
2.5	2.22044604925031e-16\\
};
\end{axis}
\end{tikzpicture}%
}
    \subfigure[\label{fig:butane:sthocastic}Stochastic {$\|P \mathbf{1} - \mathbf{1}\|_\infty$}]{\begin{tikzpicture}

\begin{axis}[%
width=0.238\columnwidth,
height=0.188\columnwidth,
at={(0\columnwidth,0\columnwidth)},
scale only axis,
xmin=0.5,
xmax=2.5,
xtick={1,2},
xticklabels={{QP-\textsc{Matlab}},{Riemannian}},
ymode=log,
ymin=1e-16,
ymax=1e-10,
yminorticks=true,
axis background/.style={fill=white},
]
\addplot [color=black, dashed, forget plot]
  table[row sep=crcr]{%
1	5.24447152372431e-12\\
1	8.77609096505694e-12\\
};
\addplot [color=black, dashed, forget plot]
  table[row sep=crcr]{%
2	9.41469124882133e-14\\
2	1.3944401189292e-13\\
};
\addplot [color=black, dashed, forget plot]
  table[row sep=crcr]{%
1	1.40387701463851e-12\\
1	2.29372076887557e-12\\
};
\addplot [color=black, dashed, forget plot]
  table[row sep=crcr]{%
2	2.70894418008538e-14\\
2	4.45199432874688e-14\\
};
\addplot [color=black, forget plot]
  table[row sep=crcr]{%
0.925	8.77609096505694e-12\\
1.075	8.77609096505694e-12\\
};
\addplot [color=black, forget plot]
  table[row sep=crcr]{%
1.925	1.3944401189292e-13\\
2.075	1.3944401189292e-13\\
};
\addplot [color=black, forget plot]
  table[row sep=crcr]{%
0.925	1.40387701463851e-12\\
1.075	1.40387701463851e-12\\
};
\addplot [color=black, forget plot]
  table[row sep=crcr]{%
1.925	2.70894418008538e-14\\
2.075	2.70894418008538e-14\\
};
\addplot [color=blue, forget plot]
  table[row sep=crcr]{%
0.85	2.29372076887557e-12\\
0.85	5.24447152372431e-12\\
1.15	5.24447152372431e-12\\
1.15	2.29372076887557e-12\\
0.85	2.29372076887557e-12\\
};
\addplot [color=blue, forget plot]
  table[row sep=crcr]{%
1.85	4.45199432874688e-14\\
1.85	9.41469124882133e-14\\
2.15	9.41469124882133e-14\\
2.15	4.45199432874688e-14\\
1.85	4.45199432874688e-14\\
};
\addplot [color=red, forget plot]
  table[row sep=crcr]{%
0.85	3.72546438143218e-12\\
1.15	3.72546438143218e-12\\
};
\addplot [color=red, forget plot]
  table[row sep=crcr]{%
1.85	6.59472476627343e-14\\
2.15	6.59472476627343e-14\\
};
\addplot [color=black, only marks, mark=+, mark options={solid, draw=red}, forget plot]
  table[row sep=crcr]{%
1	1.32174271527674e-11\\
1	1.7593260182025e-11\\
};
\addplot [color=black, only marks, mark=+, mark options={solid, draw=red}, forget plot]
  table[row sep=crcr]{%
2	1.74749104076e-13\\
2	1.96731519963578e-13\\
2	2.37365682664858e-13\\
2	3.22852855560996e-13\\
2	4.37005986952954e-12\\
2	7.65021379578457e-12\\
};
\addplot [color=black, dashed, forget plot]
  table[row sep=crcr]{%
0.5	2.22044604925031e-16\\
2.5	2.22044604925031e-16\\
};
\end{axis}
\end{tikzpicture}%
}
    
    \caption{Example of estimation of the reversible Markov chain $P$ from the $A$ obtained through the normalization of the transition count matrix for the Butane potential; the horizontal dashed lines in panels~\ref{fig:butane:detailed}, \ref{fig:butane:stationary}, and~\ref{fig:butane:sthocastic} represent the machine precision for \lstinline[language=Matlab]{double} in \textsc{Matlab}, i.e., the \lstinline[language=Matlab]{eps} value of $2^{-52} \approx 2.2 \times 10^{-16}$.}
    \label{fig:butane}
\end{figure}

The results across all evaluated metrics are summarized as \emph{box plots} in Figure~\ref{fig:butane}---the box limits indicate the range of the central $50\%$ of the data, with a central line marking the median value while the dashed whiskers extend to the most extreme data points not considered outliers which are in turn plotted individually using the \lstinline[language=Matlab]!'+'! marker symbol. As discussed in Section~\ref{sec:syntethic}, we observe that for small-sized Markov chains—corresponding in this case to the decomposition of the interval \([0, 2\pi]\) into \( \{\Omega_i\}_{i=1}^{30} \) the QP-based algorithm is faster than the Riemannian solver. {This behavior is mainly due to the relatively small problem size ($30\times30$), the overhead associated with constructing and operating with the Riemannian optimization framework represents in this case a non-negligible portion of the total runtime, whereas the quadratic programming routines such as \lstinline[language=Matlab]{quadprog} can solve problems of this scale very efficiently. For larger instances, this overhead becomes negligible relative to the cost of the iterations, and the advantages of the Riemannian approach become more apparent, as illustrated in the experiments of Section~\ref{sec:syntethic}.}

The relative distances from reversibility obtained by the two methods are nearly identical, and in both cases, the detailed balance equations are satisfied up to machine precision. Among the two approaches, Algorithm~\ref{alg:nearest_reversible} demonstrates greater numerical accuracy in satisfying both the stochasticity condition and the requirement that \( \boldsymbol{\pi} \) be the stationary distribution of the optimized transition matrix.

{{ \subsection{A problem with solution on the boundary of $\mathcal{M}_{\boldsymbol{\pi}}$}}}
\label{subsec:zero entries}

{
It is worth noting that, in some applications, the solution of problem~\eqref{eq:optimization_problem} may contain zero entries. In this setting, the Riemannian method described in Section~\ref{sec:algorithm} returns a minimizer in $\mathcal{R}_{>}$ whose sparsity pattern closely resembles that of the target matrix. As an illustrative example, we construct an $8 \times 8$ reversible matrix $A$ with a banded structure and denote by $\boldsymbol{\pi}$ its stationary distribution. This step can be obtained employing a Metropolis--Hastings adjustment~\cite{Metropolis}, to a banded stochastic matrix with stationary distribution $\boldsymbol{\pi}$.
Next, we perturb $A$ with a dense matrix $E$ of small Frobenius norm, in order to obtain a full matrix. After adding the perturbation, we renormalize the rows so that the resulting matrix $\widetilde{A}= D_{(A+E)\mathbf{1}}^{-1}(A+E)$ remains stochastic; however, this operation generally destroys reversibility. Finally, we apply the approach of Gillis and Van Dooren~\cite{Gillis}\footnote{The relevant code is available at \url{https://gitlab.com/ngillis/TSDP}.} to compute a correction matrix $\Delta$ that restores the prescribed stationary distribution, thus building a final possibly non-reversible stochastic matrix $S = \widetilde{A} + \Delta$. A possible MATLAB implementation, given $A$ and $\boldsymbol{\pi}$, is reported in the following pseudocode.
}

\begin{lstlisting}[style=Matlab-editor,basicstyle=\footnotesize,backgroundcolor=\color{MatlabCellColour}]
E = 1e-5.*rand(n); %
Atilde = A+E;
row_sums = sum(Atilde,2);
Atilde = diag(1./row_sums) * Atilde; %
[Delta,~,~,~] = supportTSDP(Atilde,pi, options); %
S = Atilde + Delta; %
\end{lstlisting}

\begin{figure}[htbp]
    \centering
    \begin{subfigure}
        \centering
        \includegraphics[width=0.48\linewidth]{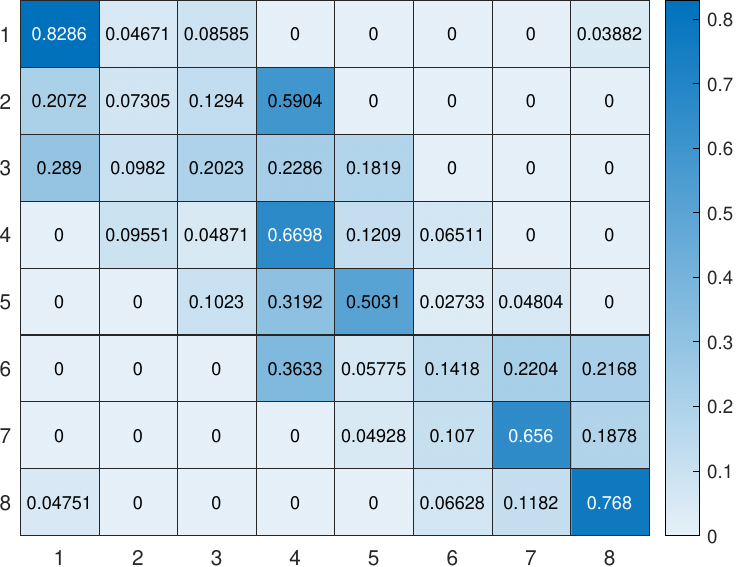}
    \end{subfigure}
    \hfill
    \begin{subfigure}
        \centering
        \includegraphics[width=0.48\linewidth]{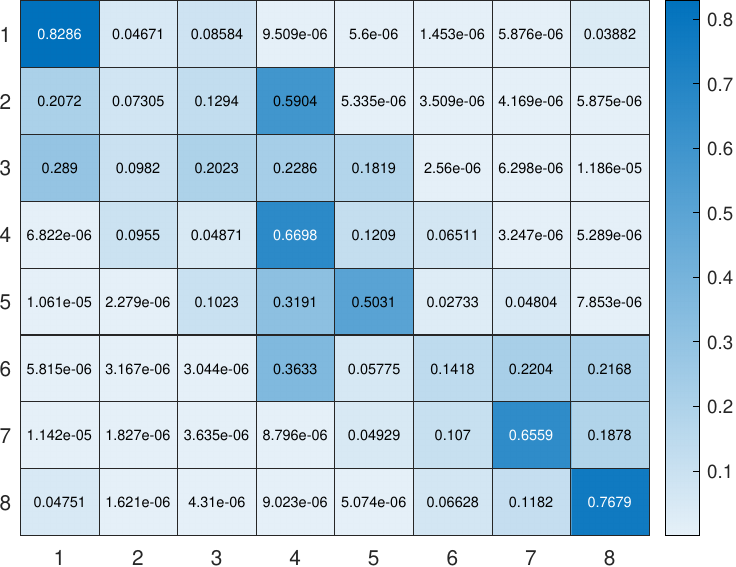}
    \end{subfigure}
    \caption{{Entries of the original reversible stochastic matrix $A$, with a banded structure and two non zero elements in $A_{18}$ and $A_{81}$ (on the left). Entries of the reversible matrix computed by the Riemannian optimizer (on the right).}}
    \label{fig:Two-Heatmap}
\end{figure}

{
We then apply the proposed method to the matrix $S$ and compute the closest reversible matrix to it. The results are shown in Figure~\ref{fig:Two-Heatmap}. Although the Riemannian method takes the full matrix $S$ as input, the output matrix $P$ recovers the banded structure of the original reversible matrix $A$. In particular, the entries outside the nonzero diagonals of the original reversible matrix $A$ have small magnitude. Finally, we report the Frobenius distances $\|P-S\|_F = 2.3971\times10^{-5}$,  $\|P-A\|_F = 6.7355\times10^{-5}$.
}

\subsection{Noisy stationary vector}

{
In this paper, given a transition matrix $A\in \mathbb{R}^{10 \times 10}$ and its stationary distribution $\boldsymbol{\pi}$, we reconstruct the reversible matrix $P$ closest to $A$, with the same stationary distribution $\boldsymbol{\pi}$. Nevertheless, in several contexts both the matrix $A$ and the stationary distribution $\boldsymbol{\pi}$ may be subjected to measurements errors and noise. In this framework, we do not have the exact stationary distribution $\boldsymbol{\pi}$, but rather its approximation $\boldsymbol{\pi}_{\text{noisy}}$. To provide an example of this instance, let us consider a reversible matrix $A$, with stationary distribution $\boldsymbol{\pi}$. Then, using $A$, we generate a discrete-time Markov chain $\left\lbrace X_t \right\rbrace_{t\ge 1}$, starting from a prescribed initial state, $X_0=1$. The simulation is performed via inverse transform sampling. From the simulated trajectory $\{X_t\}_{t=1}^{N}$, we estimate an empirical transition matrix $A_{\text{noisy}}$ for the Markov chain. This can be achieved constructing a count matrix $C$ and observing the number of transition for the state $i$ to the state $j$. The matrix $A_{\text{noisy}}$ is obtained normalizing the entries of $C$, as in Section~\ref{sec:application}, and we denote by $\boldsymbol{\pi}_{\text{noisy}}$ its stationary distribution.
}

\begin{figure}[htbp]  %
    \centering
\begin{tikzpicture}
\begin{axis}[%
width=0.9\textwidth,
height=0.3\textwidth,
scale only axis,
xmode=log,
xmin=100,
xmax=100000000,
xlabel = $N$,
ymode=log,
ymin=0.0001,
ymax=1.26798404475842,
yminorticks=true,
axis background/.style={fill=white},
legend style={legend cell align=left, align=left, draw=none, fill=none}
]
\addplot [color=mycolor1, mark=*, mark size=2pt, line width =1.5pt]
  table[row sep=crcr]{%
100	0.300787881287673\\
1000	0.158240482646136\\
10000	0.0423004101100549\\
100000	0.0121858813442963\\
1000000	0.00455495629029363\\
10000000	0.00150667926465973\\
100000000	0.000372055590797132\\
};
\addlegendentry{Obj. function with $\boldsymbol{\pi}_{\text{noisy}}$}

\addplot [color=mycolor2, mark=*, mark size=2pt, line width =1.5pt]
  table[row sep=crcr]{%
100	0.657716379510912\\
1000	0.197265252188457\\
10000	0.0468288446595037\\
100000	0.0131909603945586\\
1000000	0.00480880562735559\\
10000000	0.0016372158766128\\
100000000	0.000473575099865795\\
};
\addlegendentry{Obj. function with $\boldsymbol{\pi}$}

\addplot [color=mycolor3, mark=*, mark size=2pt, line width =1.5pt]
  table[row sep=crcr]{%
100	1.26798404475842\\
1000	0.36156222032859\\
10000	0.114273206500229\\
100000	0.0360493843029364\\
1000000	0.0109489391784444\\
10000000	0.00427153933355256\\
100000000	0.00111039211903286\\
};
\addlegendentry{Recons. error with $\boldsymbol{\pi}_{\text{noisy}}$}

\addplot [color=mycolor4, mark=*,  mark size=2pt, line width =1.5pt]
  table[row sep=crcr]{%
100	1.1390416510808\\
1000	0.342452472382728\\
10000	0.111637902000502\\
100000	0.0356149894276781\\
1000000	0.0105526983387637\\
10000000	0.00421875445214382\\
100000000	0.00106893523556125\\
};
\addlegendentry{Recons. error with $\boldsymbol{\pi}$}

\end{axis}
\end{tikzpicture}%

\caption{Behavior of the method with respect a noisy stationary distribution $\boldsymbol{\pi}_{\text{noisy}}$, with respect to the number of samples for the simulation, $N$.}\label{fig:noisy}
\end{figure}

{We run our algorithm twice: the first time employing the pair $(A_{\text{noisy}},\boldsymbol{\pi})$ and the second time with the pair $(A_{\text{noisy}},\boldsymbol{\pi}_{\text{noisy}})$. Consequently, we denote by $P$ and $P_{\text{noisy}}$ the reversible matrices, obtained by means of Algorithm~\ref{alg:nearest_reversible} and inputs $(A_{\text{noisy}},\boldsymbol{\pi})$ and $(A_{\text{noisy}},\boldsymbol{\pi}_{\text{noisy}})$, respectively.
In Figure~\ref{fig:noisy}, we report both the reconstruction errors $\|P - A\|_F$, $\|P_{\text{noisy}} - A\|_F$, and the value of the objective functional $\|P- A_{\text{noisy}}\|_F$, $\|P_{\text{noisy}}- A_{\text{noisy}}\|_F$, as a function of the sample size $N$.
}

{
The results of this example suggest that the reconstruction of the initial reversible matrix $A$ and the value of the objective functional behave in a similar way, when taking into account either the original stationary distribution or the noisy one. Observe that this holds as long as the number of considered samples $N$ is sufficiently large. Then, even in framework where we only have an estimate of the stationary distribution, the method could be employed to approximate the closest reversible matrix to a given one. Finally, for completeness, in Figure \ref{fig:lower}, we plot the lower bound previously derived in Proposition \ref{prop:lower_bound}, where we select the $\infty$-norm.
}

\begin{figure}
    \centering
 \begin{tikzpicture}

\begin{axis}[%
width=0.9\textwidth,
height=0.3\textwidth,
scale only axis,
xmode=log,
xmin=100,
xmax=100000000,
xminorticks=true,
ymode=log,
ymin=0.0001,
ymax=1,
yminorticks=true,
axis background/.style={fill=white},
legend style={legend cell align=left, align=left, draw=none, fill=none}
]
\addplot [color=mycolor1, mark=*, mark size=2pt, line width =1.5pt]
  table[row sep=crcr]{%
100	0.73311108157036\\
1000	0.11771226048448\\
10000	0.023613563834205\\
100000	0.00592115269120858\\
1000000	0.00307788777512549\\
10000000	0.000712693725878889\\
100000000	0.000485367098204507\\
};
\addlegendentry{$\| P - P_{\text{noisy}} \|_{\infty}$}

\addplot [color=mycolor2, mark=*, mark size=2pt, line width =1.5pt]
  table[row sep=crcr]{%
100	0.12496848467822\\
1000	0.061112849270684\\
10000	0.0153819653848218\\
100000	0.0020934685328534\\
1000000	0.00124040684600334\\
10000000	0.000352561028087203\\
100000000	0.000190356472601799\\
};
\addlegendentry{Lower bound w.r.t. $\| \cdot \|_\infty$}

\end{axis}
\end{tikzpicture}%

    \caption{Lower bound for $\| P - P_{\text{noisy}} \|_{\infty}$, as in Proposition~\ref{prop:lower_bound}, with respect to the number of samples $N$ for the simulation.}
    \label{fig:lower}
\end{figure}

\section{Conclusions and future perspectives}\label{sec:conclusions}

In this work, we have introduced a Riemannian optimization-based approach for computing the closest transition matrix of a reversible Markov chain to a given (possibly non-reversible) chain, with respect to the Frobenius norm. To this end, we proposed a modification of the Riemannian embedded manifold---endowed with the Fisher information metric---of stochastic matrices with a prescribed stationary distribution, by restricting it to the embedded manifold of reversible stochastic matrices. This construction can also be interpreted as a generalization of the Riemannian manifold of symmetric stochastic matrices.

The proposed numerical algorithm demonstrates robustness and efficiency when compared to existing approaches based on QP, outperforming them in terms of both memory usage and computational time.

A current limitation of the method lies in its reliance on the assumption that the output matrix is either dense or decomposable into dense ergodic classes. This is due to the fact that the Fisher information metric requires matrix elements to be strictly positive, thereby constraining the manifold to the interior of the probability simplex. {To enforce sparsity constraints, we started working\cite{cipolla2026nearestreversiblemarkovchains} on an approach based on quadratic programming for Markov chains with sparse transition matrices.} Future work will explore the relaxation of the strict positivity constraint and the solution of problems with sparse stochastic matrices, possibly incorporating the usage of lifting techniques\cite{MR4851039}. These techniques enable smooth parameterizations of sets that may be nonsmooth, potentially broadening the applicability of the proposed framework to sparse or structured stochastic matrices. Conversely, in scenarios where one can assume that the transition matrix is dense, it may be worthwhile to consider replacing the Frobenius norm in the objective function with the Kullback–Leibler (KL) divergence\cite{Wolfer2021393}. The KL divergence, like the Fisher metric, also requires strictly positive entries for its computation and naturally aligns with the information-theoretic interpretation of Markov chain estimation. 

\bmsection*{Author contributions}

All authors contributed equally.

\bmsection*{Acknowledgments}
All the authors are member of the INdAM GNCS and acknowledge the MUR Excellence Department Project awarded to the Department of Mathematics, University of Pisa, CUP I57G22000700001. The authors would like to thank the anonymous referees for the careful reading of the manuscript and for their suggestions, which have improved the quality and clarity of the paper.

\bmsection*{Financial disclosure}

None reported.

\bmsection*{Conflict of interest}

The authors declare no potential conflict of interests.

\bibliography{biblio}

\begin{thebibliography}{10}
\providecommand \doibase [0]{http://dx.doi.org/}%

\bibitem{Seneta}
Seneta E. {Markov Chains as Models in Statistical Mechanics}. {\it Statistical
  Science.} 2016\string;31(3)\string:399 -- 414.
\newblock \href {\doibase 10.1214/16-STS568} {doi: 10.1214/16-STS568}

\bibitem{PhysRevE.73.046126}
Alvarez J, Hajek B. Equivalence of trans paths in ion channels. {\it Phys. Rev.
  E.} 2006\string;73\string:046126.
\newblock \href {\doibase 10.1103/PhysRevE.73.046126} {doi:
  10.1103/PhysRevE.73.046126}

\bibitem{MR1397966}
Gilks WR, Richardson S, Spiegelhalter DJ. \kern-2pt, eds.{\it Markov chain
  {M}onte {C}arlo in practice}.
\newblock Interdisciplinary StatisticsChapman \& Hall, London, 1996

\bibitem{Wu2015184}
Wu SJ, Chu MT. Constructing optimal transition matrix for Markov chain Monte
  Carlo. {\it Linear Algebra and Its Applications.} 2015\string;487\string:184
  – 202.
\newblock \href {\doibase 10.1016/j.laa.2015.09.016} {doi:
  10.1016/j.laa.2015.09.016}

\bibitem{Fackeldey201773}
Fackeldey K, Niknejad A, Weber M. Finding metastabilities in reversible Markov
  chains based on incomplete sampling. {\it Special Matrices.}
  2017\string;5(1)\string:73 – 81.
\newblock \href {\doibase 10.1515/spma-2017-0006} {doi: 10.1515/spma-2017-0006}

\bibitem{Tifenbach2013120}
Tifenbach RM. A combinatorial approach to nearly uncoupled Markov chains I:
  Reversible Markov chains. {\it Electronic Transactions on Numerical
  Analysis.} 2013\string;40\string:120 – 147.

\bibitem{Grone2008}
Grone R, Hoffmann KH, Salamon P. An interlacing theorem for reversible Markov
  chains. {\it Journal of Physics A: Mathematical and Theoretical.}
  2008\string;41(21).
\newblock \href {\doibase 10.1088/1751-8113/41/21/212002} {doi:
  10.1088/1751-8113/41/21/212002}

\bibitem{PhysRevE.93.033307}
K\"ursten R, Behn U. Patchwork sampling of stochastic differential equations.
  {\it Phys. Rev. E.} 2016\string;93\string:033307.
\newblock \href {\doibase 10.1103/PhysRevE.93.033307} {doi:
  10.1103/PhysRevE.93.033307}

\bibitem{Choi_Wolfer}
Choi MCH, Wolfer G. Systematic Approaches to Generate Reversiblizations of
  Markov Chains. {\it IEEE Transactions on Information Theory.}
  2024\string;70(5)\string:3145-3161.
\newblock \href {\doibase 10.1109/TIT.2023.3304685} {doi:
  10.1109/TIT.2023.3304685}

\bibitem{Choi}
Choi MCH, Hird M, Wang Y. Improving the {C}onvergence of {M}arkov {C}hains via
  {P}ermutations and {P}rojections. {\it Random Structures \& Algorithms.}
  2025\string;66(4)\string:e70016.
\newblock \href {\doibase https://doi.org/10.1002/rsa.70016} {doi:
  https://doi.org/10.1002/rsa.70016}

\bibitem{Nielsen2015483}
Nielsen A, Weber M. Computing the nearest reversible Markov chain. {\it
  Numerical Linear Algebra with Applications.} 2015\string;22(3)\string:483 –
  499.
\newblock \href {\doibase 10.1002/nla.1967} {doi: 10.1002/nla.1967}

\bibitem{DurastanteMeini}
Durastante F, Meini B. Stochastic {$p$}th root approximation of a stochastic
  matrix: a {R}iemannian optimization approach. {\it SIAM J. Matrix Anal.
  Appl..} 2024\string;45(2)\string:875--904.
\newblock \href {\doibase 10.1137/23M1589463} {doi: 10.1137/23M1589463}

\bibitem{Douik8861409}
Douik A, Hassibi B. Manifold Optimization Over the Set of Doubly Stochastic
  Matrices: A Second-Order Geometry. {\it IEEE Transactions on Signal
  Processing.} 2019\string;67(22)\string:5761-5774.
\newblock \href {\doibase 10.1109/TSP.2019.2946024} {doi:
  10.1109/TSP.2019.2946024}

\bibitem{KemenySnell}
Kemeny JG, Snell JL. {\it Finite {M}arkov chains}.
\newblock The University Series in Undergraduate MathematicsD. Van Nostrand
  Co., Inc., Princeton, N.J.-Toronto-London-New York, 1960.

\bibitem{NocedalWright}
Nocedal J, Wright SJ. {\it Numerical optimization}.
\newblock Springer Series in Operations ResearchSpringer-Verlag, New York, 1999

\bibitem{MR4533407}
Boumal N. {\it An introduction to optimization on smooth manifolds}.
\newblock Cambridge University Press, Cambridge, 2023.

\bibitem{AbsilBook}
Absil PA, Mahony R, Sepulchre R. {\it Optimization algorithms on matrix
  manifolds}.
\newblock Princeton University Press, Princeton, NJ, 2008.
\newblock With a foreword by Paul Van Dooren

\bibitem{Lee2018}
Lee JM. {\it Introduction to {R}iemannian Manifolds}.
\newblock Springer International Publishing, 2018

\bibitem{MR304178}
Tarjan R. Depth-first search and linear graph algorithms. {\it SIAM J.
  Comput..} 1972\string;1(2)\string:146--160.
\newblock \href {\doibase 10.1137/0201010} {doi: 10.1137/0201010}

\bibitem{MR1800071}
Amari Si, Nagaoka H. {\it Methods of information geometry}. 191 of {\it
  Translations of Mathematical Monographs}.
\newblock American Mathematical Society, Providence, RI; Oxford University
  Press, Oxford, 2000.
\newblock Translated from the 1993 Japanese original by Daishi Harada

\bibitem{MR544666}
Berman A, Plemmons RJ. {\it Nonnegative matrices in the mathematical sciences}.
  9 of {\it Classics in Applied Mathematics}.
\newblock Society for Industrial and Applied Mathematics (SIAM), Philadelphia,
  PA, 1994.
\newblock Revised reprint of the 1979 original

\bibitem{GolubAndVanLoan}
Golub GH, Van~Loan CF. {\it Matrix computations}.
\newblock Johns Hopkins Studies in the Mathematical SciencesJohns Hopkins
  University Press, Baltimore, MD.
\newblock fourth~ed., 2013.

\bibitem{Sinkhorn}
Sinkhorn R. {A Relationship Between Arbitrary Positive Matrices and Doubly
  Stochastic Matrices}. {\it The Annals of Mathematical Statistics.}
  1964\string;35(2)\string:876 -- 879.
\newblock \href {\doibase 10.1214/aoms/1177703591} {doi:
  10.1214/aoms/1177703591}

\bibitem{MR210731}
Sinkhorn R, Knopp P. Concerning nonnegative matrices and doubly stochastic
  matrices. {\it Pacific J. Math..} 1967\string;21\string:343--348.

\bibitem{MR179182}
Marcus M, Newman M. Generalized functions of symmetric matrices. {\it Proc.
  Amer. Math. Soc..} 1965\string;16\string:826--830.
\newblock \href {\doibase 10.2307/2033933} {doi: 10.2307/2033933}

\bibitem{MR2399579}
Knight PA. The {S}inkhorn-{K}nopp algorithm: convergence and applications. {\it
  SIAM J. Matrix Anal. Appl..} 2008\string;30(1)\string:261--275.
\newblock \href {\doibase 10.1137/060659624} {doi: 10.1137/060659624}

\bibitem{manopt}
Boumal N, Mishra B, Absil PA, Sepulchre R. {M}anopt, a {M}atlab Toolbox for
  Optimization on Manifolds. {\it Journal of Machine Learning Research.}
  2014\string;15(42)\string:1455--1459.

\bibitem{MR2335248}
Absil PA, Baker CG, Gallivan KA. Trust-region methods on {R}iemannian
  manifolds. {\it Found. Comput. Math..} 2007\string;7(3)\string:303--330.
\newblock \href {\doibase 10.1007/s10208-005-0179-9} {doi:
  10.1007/s10208-005-0179-9}

\bibitem{nr}
Rossi RA, Ahmed NK. {The Network Data Repository with Interactive Graph
  Analytics and Visualization}. In:  2015.

\bibitem{davis1941deep}
Davis A, Gardner BB, Gardner MR. Deep {S}outh {C}hicago.;  1941.

\bibitem{SciPyProceedings_11}
Hagberg AA, Schult DA, Swart PJ. {Exploring Network Structure, Dynamics, and
  Function using NetworkX}. In:  Varoquaux G, Vaught T, Millman J. \kern-2pt,
  eds. {\it Proceedings of the 7th Python in Science Conference}\,SciPy 2008.
  2008; Pasadena, CA USA\string:11 - 15.

\bibitem{MR1875515}
Dolan ED, Mor\'e JJ. Benchmarking optimization software with performance
  profiles. {\it Math. Program..} 2002\string;91(2)\string:201--213.
\newblock \href {\doibase 10.1007/s101070100263} {doi: 10.1007/s101070100263}

\bibitem{Kloeden1992}
Kloeden PE, Platen E. {\it Stochastic Differential Equations}\string:103--160;
  Berlin, Heidelberg: Springer Berlin Heidelberg .
\newblock 1992

\bibitem{Metropolis}
Metropolis N, Rosenbluth AW, Rosenbluth MN, Teller AH, Teller E. {Equation of
  State Calculations by Fast Computing Machines}. {\it J. Chem. Phys..}
  1953\string;21(6)\string:1087-1092.
\newblock \href {\doibase 10.1063/1.1699114} {doi: 10.1063/1.1699114}

\bibitem{Gillis}
Gillis N, Van~Dooren P. Assigning stationary distributions to sparse stochastic
  matrices. {\it SIAM J. Matrix Anal. Appl..}
  2024\string;45(4)\string:2184--2210.
\newblock \href {\doibase 10.1137/23M1627328} {doi: 10.1137/23M1627328}

\bibitem{cipolla2026nearestreversiblemarkovchains}
Cipolla S, Durastante F, Gnazzo M, Meini B. Nearest Reversible {M}arkov Chains
  with Sparsity Constraints: An Optimization Approach.;  2026.

\bibitem{MR4851039}
Levin E, Kileel J, Boumal N. The effect of smooth parametrizations on nonconvex
  optimization landscapes. {\it Math. Program..}
  2025\string;209(1-2)\string:63--111.
\newblock \href {\doibase 10.1007/s10107-024-02058-3} {doi:
  10.1007/s10107-024-02058-3}

\bibitem{Wolfer2021393}
Wolfer G, Watanabe S. Information Geometry of Reversible Markov Chains. {\it
  Information Geometry.} 2021\string;4(2)\string:393 – 433.
\newblock \href {\doibase 10.1007/s41884-021-00061-7} {doi:
  10.1007/s41884-021-00061-7}

\end{thebibliography}

\end{document}